\documentclass[12pt]{article}
\usepackage{amsmath,amssymb,euscript}
\def\sect#1{\advance\sectionnum by 1\section{#1}}
\newcount\sectionnum \sectionnum=-1
\def\thesection{\the\sectionnum}

\newtheorem{theorem}{Theorem}[section]
\newtheorem{lemma}[theorem]{Lemma}
\newtheorem{definition}[theorem]{Definition}
\newtheorem{proposition}[theorem]{Proposition}
\newtheorem{corollary}[theorem]{Corollary}
\newtheorem{remark}[theorem]{Remark}

\newcommand{\eb}{\begin{equation}}
\newcommand{\ee}{\end{equation}}

\begin{document}

\title{\bf  Topologies on the group of homeomorphisms of a Cantor set}
\author{{\bf S.~Bezuglyi}\thanks{Supported in part by  CRDF grant
UM1-2546-KH-03}\\ Institute for Low Temperature Physics, Kharkov, Ukraine\\
\and{\bf A.H.~Dooley}\\ University of New South Wales, Sydney, Australia\\
\and {\bf J.~Kwiatkowski\thanks{The research was supported by the grant KBN
5 P03A 027 21}}\\ Nicolas Copernicus University, Torun, Poland}

\date{}

\maketitle

\begin{abstract}
Let $Homeo(\Omega)$ be the group of all homeomorphisms of a Cantor set
$\Omega$. We study topological properties of $Homeo(\Omega)$ and its subsets
with respect to the uniform $(\tau)$ and weak $(\tau_w)$ topologies. The
classes of odometers and periodic, aperiodic, minimal, rank 1 homeomorphisms
are considered and the closures of those classes in $\tau$ and $\tau_w$ are
found.
\end{abstract}

\newcommand{\La}{\Lambda}
\newcommand{\Om}{\Omega}
\newcommand{\om}{\omega}
\newcommand{\e}{\varepsilon}
\newcommand{\al}{\alpha}
\newcommand{\vp}{\varphi}
\newcommand{\G}{\Gamma}
\newcommand{\g}{\gamma}
\newcommand{\N}{{\mathbb N}}
\newcommand{\T}{\theta}
\newcommand{\De}{\Delta}
\newcommand{\de}{\delta}
\newcommand{\s}{\sigma}
\newcommand{\E}{{\cal E}}
\newcommand{\A}{{\cal A}}
\newcommand{\F}{{\cal F}}
\newcommand{\B}{{\cal B}}
\newcommand{\M}{{\cal M}}
\newcommand{\R}{{\cal R}}
\newcommand{\bs}{(X,{\cal B})}
\newcommand{\Aut}{Aut(X,{\cal B})}
\newcommand{\h}{Homeo(\Om)}
\newcommand{\pos}{{\mathbb R}^*_+}
\newcommand{\la}{\lambda}
\newcommand{\Z}{{\mathbb Z}}
\newcommand{\ap}{{\cal A}p}
\newcommand{\per}{{\cal P}er}
\newcommand{\od}{{\cal O}d}
\newcommand{\Min}{{\cal M}in}
\newcommand{\mix}{{\cal M}ix}
\newcommand{\mov}{{\cal M}ov}
\newcommand{\wm}{w{\cal M}ov}

\sect{Introduction}

The present paper is a continuation of our
article [BDK] about topologies on the group $\Aut$ of all Borel
automorphisms of a standard Borel space. In the introduction to that
article, we discussed our approach to the study  of topologies on
groups of transformations of an underlying space. As we mentioned
there, we were motivated, first of all, by remarkable results in
ergodic theory concerning topological properties of the group of all
automorphisms of a standard measure space. We refer to the classical
articles of Halmos [H] and Rokhlin [R] where the uniform and weak
topologies appeared as \textquotedblleft key players\textquotedblright in
ergodic theory.

The central object of the present paper is the group $\h$ of all
homeomorphisms of a Cantor set $\Om$. Although we consider several
topologies on $\h$, this group is mostly studied under two topologies,
$\tau$ and $\tau_w$. These are analogues of the uniform and weak topologies
in measurable dynamics. We should remark that $\tau_w$ is, in fact, the
usual sup-topology of uniform convergence which has occurred in many papers
on topological dynamics (see, e.g. [GPS2, GW2]). Many interesting questions
can be asked about the topological properties of $\h$ and its subsets. For
instance, E.~Glasner and B.~Weiss studied in [GW2] the Rokhlin property on
$(\h, \tau_w)$ showing that the action of $\h$ on itself by conjugation has
dense orbits. In this article, we will concentrate on the following
directions, which we believe are natural initial questions in this theory:
(i) global properties of some basic topologies on $\h$, (ii) finding
closures of subsets of $\h$ consisting of periodic, aperiodic, minimal,
topologically transitive, rank 1 homeomorphisms, and odometers in $\tau$ and
$\tau_w$.

It might be asked why we consider only  Cantor sets as the underlying
space.  First of all, we remark that Cantor sets and their
homeomorphisms arise naturally in various areas of dynamical systems, for
example in fractals, low-dimensional dynamics etc. Although  topological
and measurable dynamics are, strictly speaking, completely different
theories, we believe that Cantor dynamics has several features in common
with measurable dynamics. To support this point of view we refer to the
results on orbit equivalence of minimal homeomorphisms and full groups
proved in [GPS1, GPS2, GW1, HPS, BK1, BK2]. We believe that the
following properties of Cantor sets underlie this similarity: (a) all
Cantor sets are homeomorphic; (b) for every Cantor set, there exists a
countable family of clopen sets generating the topology; (c) any Cantor set
can be partitioned into a finite collection of clopen subsets.
Nevertheless, we are optimistic that some ideas of this paper may be used
in the context of general topological dynamics.

The paper is organized as follows. In Section 1, we introduce several
topologies on $\h$ and study global topological properties of $\h$ mostly
with respect to $\tau$ and $\tau_w$. All possible relations between these
topologies are found in Theorem 1.3. We mention the curious fact that
$(\h,\tau_w)$ is a zero-dimensional Polish space. It turns out that $\tau_w$
is equivalent to the topology $p$ whose base of neighborhoods is defined by
$W(T; F_1,...,F_k) = \{S\in \h\ |\ SF_i = TF_i, \ i= 1,...,k\}$ where $F_i$
is clopen. This fact is a justification of the name \textquotedblleft weak"
topology which we use for $\tau_w$. Section 2 deals principally with the
problem of approximation by periodic homeomorphisms. We prove a topological
version of the Rokhlin lemma for minimal homeomorphisms for both $\tau$ and
$\tau_w$. On the other hand, we show that pointwise periodic homeomorphisms
are not dense in $(\h, \tau_w)$. Amongst other results, we obtain a
description of periodic and aperiodic homeomorphisms from the topological
full group of a minimal homeomorphism. In Section 3, we consider
homeomorphisms of rank 1 and show that they are necessarily odometers. In
the last section, we study closures of various subsets in $\h$ with respect
to $\tau$ and $\tau_w$. In particular, we prove that the $\tau_w$-closure of
the set of minimal homeomorphisms is the same as the closure of the set of
odometers. Moreover, we give a dynamical description of homeomorphisms which
belong to the closure: $T\in \h$ belongs to the $\tau_w$-closure of the set
of minimal homeomorphisms if and only if $T$ has the following property: for
every non-trivial clopen $F$, the sets $TF\setminus F$ and $F\setminus TF$
are non-empty.
\smallskip

Throughout the paper, we use the following standard {\bf notation}:

\begin{itemize}

\item $\Om$ is a Cantor set;

\item $CO(\Om)$ is the family of all clopen subsets in $\Om$;

\item $\h$ is the group of all homeomorphisms of $\Om$ with  identity map
${\mathbb I}\in \h$;

\item $Aut(X, \B)$ is the group of all one-to-one Borel automorphisms of a
standard Borel space $(X,\B)$;

\item $\ap$ is the set of all aperiodic homeomorphisms;

\item $\per$ is the set of all pointwise periodic homeomorphisms and
$\per_0$ is the  subset of $\per$ consisting of  homeomorphisms with
finite period;

\item $\Min$ is the set of all minimal homeomorphisms;

\item $\mix$ is the set of all mixing homeomorphisms;

\item ${\cal M}_1(\Om)$ is the set of all Borel probability measures on
$\Om$;

\item $\de_x$ is the Dirac measure at $x \in \Om$;

\item $E(S,T)=\{ x\in \Om\ |\ Tx\ne Sx\} \cup \{x\in X \ |\ T^{-1}x\ne
S^{-1}x\}$ where $S,T \in \h$;

\item $\mu(f)=\int_X f\,d\mu$ where $f$ is in $C(\Om)_1$ (= the set of
continuous real-valued functions with $\| f \| :=\sup\{|f(x)| \ :\ x\in
\Om\} \le 1$), and $\mu\in {\cal M}_1(\Om)$;

\item $\mu\circ S(A):=\mu (SA) $ and $\mu \circ S(f):=\int_{\Om}f\,d(\mu
\circ S) = \int_{\Om}f(S^{-1}x)\,d\mu (x)$ where $S\in \h$;

\item $A^c = \Om \setminus A$.

\end{itemize}

\sect{Topologies on $Homeo(\Om)$}

\setcounter{equation}{0}

In this section, we define several topologies on $\h$. These topologies
are similar to those studied in [BDK] for $\Aut$. We make the
following (rather obvious) changes to the settings of [BDK]:   a
standard Borel space $(X,\B)$ is replaced by a Cantor set $\Om$, and Borel
sets and functions are replaced by clopen sets and continuous functions.

\begin{definition}\label{def}  {\rm (cf. [BDK, Definition 1.1])}
(i) The uniform topology $\tau$ on $\h$ is defined as
the relative topology on $\h$ induced from $(Aut(\Om,\B), \tau)$. The base
of neighborhoods is formed by
\begin{equation}\label{tau}
U(T; \mu_1,...,\mu_n; \e)= \{S\in \h)\ |\ \mu_i(E(S,T))<\e ,\ i=1,...,n\}.
\end{equation}
(ii) The topology $\tau'$ is defined on $\h$ by the base of neighborhoods
\begin{equation}\label{tau'}
U'(T; \mu_1,...,\mu_n; \e )= \{ S\in \h \ |\ \sup_{F\in CO(\Om)}\mu_i(TF\
\De \ SF)<\e,\ i=1,...,n\}.
\end{equation}
(iii) The topology $\tau''$ is defined on $\h$ by the base of neighborhoods
\begin{equation}\label{tau''}
\begin{array}{ll}
U''(T; \mu_1,...,\mu_n; \e ) = & \{ S\in \h \ |\ \sup_{f\in C(\Om)_1} \vert
\mu_i\circ S(f)-\mu_i\circ T(f)\vert <\e,\\ & i=1,...,n\}.
\end{array}
\end{equation}
(iv) The topology $p$ is defined on $\h$ by the base of neighborhoods
\begin{equation}\label{p}
W(T; F_1,...,F_k) = \{S\in \h\ |\ SF_i = TF_i, \ i= 1,...,k\}.
\end{equation}
(v) The topology $\bar p$ is defined on $\h$ by the base of neighborhoods
\begin{equation}\label{barp}
\begin{array}{ll}
&\overline{W}(T; F_1,...,F_k; \mu_1,...,\mu_n; \e) = \{S\in \h\ |\\
&\mu_j(SF_i\ \De\ TF_i) + \mu_j(S^{-1}F_i\ \De\ T^{-1}F_i) < \e,\ i=1,...,n;
 j=1,...,k\}.
\end{array}
\end{equation}
In all the above definitions, we have taken  $T\in \h$, $\mu_i \in {\cal
M}_1(\Om)$, and $F_i\in CO(\Om),\ i=1,...,n$.
\end{definition}

It is a simple exercise to verify that the collections of sets
so defined do indeed form bases of topologies.

We remark that in the Borel case (see [BDK]) we have also defined the
topologies $\tau_B, \tau'_B, \tau''_B, p_B$, and $\bar p_B$ (the subindex
$B$ stands for Borel; in [BDK] these topologies were denoted without $B$).
In fact, only one of them, $\tau$, is the relative topology induced on $\h$
from $(Aut(\Om, \B), \tau_B)$. The others are not relative topologies on
$\h$ because in their definition we use clopen subsets and continuous
functions instead of Borel ones (see (\ref{tau'}) - (\ref{barp})).
Obviously, $\tau'_B, \tau''_B, p_B$, and $\bar p_B$ being induced onto $\h$
from $Aut(\Om,\B)$ are at least not weaker than the corresponding topologies
$\tau', \tau'', p$, and $\bar p$. Thus, we have to deal here with the
topological counterparts of topologies studied in [BDK]. Nevertheless, we
will see that the greater part of our results about relations between the
topologies proved in [BDK] is still true in the context of homeomorphisms of
Cantor sets. In most cases, the proofs for homeomorphisms are either
word-for-word repetition of those in the Borel case or can be easily
adapted.

\begin{definition}\label{weaktop} For $S,T\in Homeo(\Om)$, define
\begin{equation}\label{weak}
d_w(S,T)=\sup_{x \in \Om}d(Sx, Tx)+\sup_{x\in \Om}d(S^{-1}x,T^{-1}x).
\end{equation}
Denote by $\tau_w$ the topology on $\h$ generated by the metric $d_w$.
\end{definition}

The topology $\tau_w$ is well known in topological dynamics and probably is
generally considered as the most natural topology on $\h$. In particular, it
can easily be seen that $\tau_w$ is equivalent to the topology defined by
the base of neighborhoods $\widetilde W(T;f_1,...,f_n;\e) = \{S\in \h :
\Vert f_i\circ T - f_i\circ S\Vert <\e,\ i=1,...,n\}$ where $f_1,...,f_n$
are $\Z$-valued continuous functions on $\Om$. The proof of this fact is
similar to that of [Theorem 3.6, BDK] for Borel dynamics.

We call it the {\it weak topology} following our point of view explained in
[BK1, BK2] (see also Theorem \ref{top-eq} (4) below).

It is well known that $(Homeo(\Om),\tau_w)$ is a Polish space (for every
compact metric space $\Om$). By $B_{\de }({T})$, we denote the set
$\{S\in \h\ | \ d_w(S,T)<\de\},\ T\in \h$.

Our first main result is the following

\begin{theorem}\label{top-eq}
$(1)$ The topologies $\tau$ and $\tau'$ are equivalent.\\ $(2)$ The
topology $\tau\ (\sim \tau')$ is strictly stronger than $\tau''$.\\ $(3)$
The topology $\tau\ (\sim \tau')$ is strictly stronger than $\bar p$.\\
$(4)$ The topology $\tau_w$ is equivalent to $p$.\\ $(5)$ The topology $p\
(\sim \tau_w)$ is strictly stronger than $\bar p$.\\ $(6)$ The topology
$\tau$ is not comparable with $\tau_w\ (\sim p)$ and the topology $\tau''$
is not comparable with $\bar p$. \end{theorem}

\noindent{\it Proof}. A direct analogue of this theorem was proved in the
context of Borel dynamics in [BDK]. The principal difference is here that
one needs to work with clopen sets instead of Borel sets. We will indicate
only what modifications need to be made for the use of $\h$.

(1) We follow the idea of the proof of [Theorem 3.1, BDK]. The fact that
$\tau \succ \tau'$ may be proved as in [BDK]. To prove that $\tau' \succ
\tau$, we show that each neighborhood $U = U({\mathbb I};
\mu_1,...,\mu_n;\e)$ contains  $ U' = U'({\mathbb I};
\mu_1,...,\mu_n;\e/2)$. By definition, $T \in  U'({\mathbb I};
\mu_1,...,\mu_n;\e/2)$ if $\mu_i(TF\ \De\ F) < \e/2$ for all clopen $F$ and
all $i= 1,...,n$. We first note that $E(T, {\mathbb I})$ is open,
$T$-invariant, and contains some clopen $E_0$ such that $E_0\cap TE_0
=\emptyset$. Thus one of the following alternatives must hold:

(a) for every clopen $F\subset  E(T, {\mathbb I})$ with $F\cap
TF=\emptyset$ there exists a clopen set $F'\supset F$ such that $F'\cap
TF'=\emptyset$; or

(b) there exists a clopen $F_0 \subset  E(T, {\mathbb I}),\ F_0\cap TF_0=
\emptyset$ which cannot be extended to a large set $F'$ preserving
disjointness of $F'$ and $TF'$ (in other words, each clopen set
$F' \supset F_0$ has a nonempty intersection with $TF')$.

Clearly, condition (b) is equivalent to the following property: $F_0\cup
TF_0\cup T^2F_0 = E(T, {\mathbb I})$. Therefore, in this case, $\mu_i(E(T,
{\mathbb I})) \le  \mu_i(TF_0\ \De\ F_0) + \mu_i(TF_0\ \De\ T^2F_0) < \e$.

If (a) holds, then one can find $F_1\in CO(\Om) $ such that $F_1\cap TF_1
=\emptyset$ and $\mu_i(E(T, {\mathbb I}) - (F_1\cup TF_1)) < \e/2$. This
implies that  $\mu_i(E(T, {\mathbb I})) < \e,\ i=1,...,n$, and by
(\ref{tau}), we are done.

(2) The method of proof that $\tau$ is strictly stronger than
$\tau''$ is the same as in [Theorem 3.5, BDK]. To show that a statement
analogous to Proposition 3.4 of [BDK] holds, we need to use clopen sets in
the definitions of two auxiliary topologies $\tilde\tau$ and $\bar\tau$ as
well as continuous functions in the proof of that proposition
(see Remark 3.7 [BDK]).  In particular, one sees that the topology
$\tau''$ is equivalent on $\h$ to the topology $\bar\tau$ defined by the
base
\begin{equation}\label{bartau}
\overline V(T; \mu_1,...,\mu_n; \e )
=  \{S\in \Aut\ \vert \sup_{F\in {CO(\Om)}}\vert\mu_j(TF) -
\mu_j(SF)\vert <\e,\ j=1,...,n\},
\end{equation}
where $T\in \h$ and $\mu_i\in  {\cal M}_1(\Om)$.

(3) The proof is a word for word repetition of [Proposition 3.2, BDK].

(4) Fix some $\de >0$ and let ${\cal  Q}=(F_i)_{i=1}^n$ be a partition
of $\Om$ into clopen sets such that ${\rm diam}(F_i)<\de,\ i=1,...,n$. If
$S\in W({\mathbb I};\ F_1,...,F_n)$, then $SF_i=F_i$, and therefore
$\sup_{x\in \Om} d(Sx,x)+\sup_{x\in \Om} d(S^{-1}x,x)\le 2\de$. This
proves that $B_{2\de}({\mathbb I})\supset W({\mathbb I};F_1,...,F_n)$.

Conversely, let $W({\mathbb I};\ F_1,...,F_n)$ be given. Take the
partition ${\cal Q}=(E_i)_{i\in I}$ which is generated by all $F_i$ and
$F_i^c=\Om-F_i,\ i=1,...,n$. Take $\e >0$ such that
 $$
 \e <\min \{ \min_{i\ne j} {\rm dist}(E_i,E_j),\ \min_i ({\rm diam}
(E_i))\}.
 $$
 Then, every $S\in B_{\e }({\mathbb I})$ has the property $SE_i=E_i$,
i.e. ${\cal Q}$ is fixed. Therefore $SF_k=F_k,\ k=1,...,n$, because every
$F_k$ is a union of some $E_i$'s. Thus, $S\in W({\mathbb I};\ F_1,...,F_k)$.

(5) As an immediate corollary of the equivalence proved in (4), we obtain
$\tau_w \succ \bar p$. To see that $\bar p$ is strictly weaker than
$\tau_w$, we note that $\bar p$ is weaker than the topology $\tau' \sim
\tau$.  If we assumed that $\bar p$ was equivalent to $\tau_w$, we would
have that $\tau$ is always stronger than $\tau_w$. But the latter is false
(see (6) or [BK1]).

(6) See [Proposition 3.6, BDK] where the pairs $\tau''$ and $\bar p$ have
been considered. The fact that $\tau$ and $\tau_w$ are not comparable
is a direct consequence of [Theorem 4.8, BK1].  \hfill${\square}$
\\

Now we formulate several statements concerning topological properties of
$\h$.

\begin{proposition}\label{property} $(1)$ $(\h,p)$ is a $0$-dimensional
complete metric space with no isolated points.\\ $(2)$ $\h$ is a Hausdorff
topological group with respect to the topologies $\tau,\tau',\tau'', p,
\tau_w$.\\ $(3)$ $\h$ is not closed in $(Aut(\Om,\B), \tau)$.\\ $(4)$ $\h$
is dense in $(Aut(\Om,\B), \tau)$.
\end{proposition}

\noindent
{\it Proof}. The proof of the first statement  follows easily
from [Proposition 1.6, BDK],  replacing Borel sets by clopen sets. In
fact it  can be shown that {\it the sets $W(T;F_1,...,F_n)$ {\rm (see
(\ref{p}))} are closed with respect to the topologies $\tau, \tau'', p$
and $\bar p$}. The second statement of the proposition is based on a
routine verification (see [BDK]). The third assertion is taken from [BK1].

(4) We need to show that for any Borel automorphism $T$ of $(\Om, \B)$, for
any $\e >0$, and for any $\mu_1,...,\mu_n \in \M_1(\Om)$ there exists a
homeomorphism $S$ of $\Om$ such that $\mu_i(E(S,T)) <\e,\ 1=1,...,n$. By
Lusin's theorem, we can find a closed subset $F_i$ of $\Om$ such that the
restriction of $T$ to $F_i$ is a one-to-one continuous map from $F_i$ onto
$T(F_i)$ and
$$
\mu_i(\Om \setminus F_i)< \frac{\e}{2},\ \  \mu_i\circ T(\Om \setminus F_i)
< \frac{\e}{2}, \ \ i =1,...,n.
$$
Let $F = \bigcup_{i=1}^n F_i$. Then $F$ is closed, $T$ is continuous on
$F$, and $\mu_i(\Om \setminus F) + \mu_i(\Om \setminus T(F)) < \e$ for all
$i$.

Since $F$ and $TF$ are closed, we can represent $\Om \setminus F$ and  $\Om
\setminus TF$ as unions of infinitely many clopen sets:  $\Om \setminus F =
\bigcup_{j=1}^\infty A_j$ and $\Om \setminus TF = \bigcup_{j=1}^\infty
A'_j$. Then by Theorem 1 of [KR]\footnote[1], the continuous map $T : F\to
TF$ can be extended to a homeomorphism $S$ of $\Om$ such that $Tx =Sx,\
x\in F$ and $T^{-1}x = S^{-1}x,\ x\in TF$. Clearly, $\mu_i(E(S,T)) < \e,\
i=1,...,n$. \hfill$\square$
\\

\noindent {\bf Convention} As mentioned above, $\h$ is not closed in
$Aut(\Om,\B)$ in the uniform topology $\tau$, therefore the $\tau$-closure
of a subset $Y\subset \h$ does not belong to $\h$, in general. For
convenience, we will  use the following convention $\overline{Y}^\tau
:=\overline{Y}^\tau \cap \h$ without further explanation.

\begin{proposition}\label{converg}  Let $(T_n)$ be a sequence of
homeomorphisms of $\Om$. Then:
\\

$(1)$ $T_n\stackrel{\tau}{\longrightarrow}S \iff \forall x \in \Om\ \exists
n(x)\in \N$ such that $\forall n >n(x), \ T_nx = Sx$.
\\

$(2)$ $T_n\stackrel{p} {\longrightarrow} {\mathbb I} \iff \forall F
\in CO(\Om) \ \exists n(F)$ such that $\forall n > n(F),\ T_nF=F$;
\\

$(3)$ $T_n\stackrel{\overline p} {\longrightarrow} {\mathbb I}$ if
and only if  $\forall \mu \in {\cal M}_1(\Om) \ \forall F \in CO(\Om)$
\begin{equation}\label{conv-in-meas}
 \mu(T_nF\ \De\ F)+\mu (T^{-1}_nF\ \De\ F)\to 0
\end{equation}
or if and only if $\forall F \in CO(\Om)$,
\begin{equation}\label{limsup}
  F=\limsup_{n\to\infty} T_nF =   \limsup_{n\to\infty} T_n^{-1}F,
\end{equation}
where
$$
\limsup_{n\to\infty} T_nF = \bigcup_m \bigcap_{n>m}T_nF.
$$
\end{proposition}

\noindent{\it Proof}. Notice that (1) is proved in [BK1] and  (2) is
obvious. Relation (\ref{conv-in-meas}) is a direct consequence of the
definitions.   To prove the other equivalence in (3), we note that for any
$x\in \Om$ and $F\in CO(\Om)$, the convergence $T_n\stackrel{\overline p}
{\longrightarrow} {\mathbb I}$ implies that $$
\de_x (T_nF\ \De \ F)+\de_x (T^{-1}F\ \De\ F) \to 0
$$
as $n\to \infty$. This means that if $x\in F$, then there exists
$n_0=n_0(x,F)$ such that $x\in T_nF$ and $x\in T^{-1}_nF$ for all $n>n_0$.
Thus, we have proved that $F \subset \bigcup_m \bigcap_{n>m} T_nF$ and
$F\subset \bigcup_m \bigcap_{n>m} T^{-1}_nF$. In fact, these inclusions
are equalities. Indeed, if we assume that there exists $x_0\in F^c=\Om-F$
with $x_0\in \bigcap_{n>m}T_nF$ for some $m$, then we have a contradiction
to the fact that $x_0$ also belongs to $\bigcup_k \bigcap_{n>k}T_nF^c$.
Thus, (\ref{limsup}) holds.

Conversely, let $E_m=\bigcap_{n>m}T_nF$ and $\bigcup_m E_m=F$. Since
$E_m\subset E_{m+1}$, we see that for any measure $\mu \in {\cal
M}_1(\Om),\ \mu E_m\to \mu F\ (m\to \infty)$. Remark that $E_m\subset
T_nF$ for all $n>m$. Therefore $E_m=E_m\cap T_nF\subset F\cap T_nF\subset
F$. Thus we have $\mu (F\cap T_nF)\to \mu F$ as $n\to \infty$. Similarly
$\mu (F\cap T^{-1}_nF)\to \mu F$. By (\ref{conv-in-meas}), the proof is
complete. \hfill{$\square$}

%%%%%%%%%%%%%%%%%%%%%%%%%%%%SECTION %%%%%%%%%%%%%%%%%%

\sect {Periodic approximation}

%%%%%%%%%%%%%%%%%%%%%%%%%%%%%%%%%%%%%%

\setcounter{equation}{0}

Let $\Om$ be a Cantor set equipped with a metric $d$ compatible with the
clopen topology. It is natural to distinguish two principal classes of
homeomorphisms of $\Om$, the periodic and the aperiodic. We will say  that
$P\in \h$ is {\it pointwise periodic} if every $P$-orbit is finite. If
$T\in \h$ has no periodic points, then $T$ is called {\it aperiodic}.
Denote these classes by $\per$ and $\ap$, respectively.

In the paper [GW2] a new interesting notion  of simple homeomorphisms was
defined. Recall that,  by definition, $S\in \h$ is {\it simple} if it
satisfies the following conditions.

(i) There exist clopen subsets $F_j$ and integers $r_j\ge 1,\ j=1,...,k,$
such that the collection  $\{S^iF_j : i=0,1,...,r_j,\ j=1,...,k\}$ is
pairwise disjoint and $S$ has period $r_j$ on $F_j$.

(ii) There exist clopen subsets $C_s,\ s=1,...,l,$ and, for each $s$, two
disjoint periodic orbits $(y_s^+, Sy_s^+,...,S^{q^+_s-1}y_s^+)$ and
$(y_s^-, Sy_s^-,...,S^{q^-_s-1}y_s^-)$ such that the sets $(S^nC_s : n\in
\Z, s=1,...,l)$ are pairwise disjoint and spiral towards the periodic
orbits of $y_s^+$ and  $y_s^-$, that is $\lim_{n\to \pm \infty}{\rm
dist}(S^nC_s, S^ny_s^{\pm}) = 0$

(iii) The space $\Om$ may be represented as
\begin{equation}\label{simple}
\Om = \bigcup_{j=1}^k\bigcup_{i=0}^{r_j-1} S^iF_j\cup \
\bigcup_{s=1}^l\bigcup_{n\in \Z} S^nC_s\ \cup\ \bigcup_{s=1}^l
\left[(y_s^+,...,S^{q^+_s-1}y_s^+) \cup (y_s^-,
...,S^{q^-_s-1}y_s^-)\right]
 \end{equation}

It was shown in [Theorem 2.2, GW2] that the set ${\cal S}$ of simple
homeomorphisms is dense in $(\h,\tau_w)$.

\begin{theorem}\label{ap}  $(1)$  $\ap$  is closed in $\h$ with respect
to the topologies $\tau$ and $\tau''$.\\
$(2)$  $\ap$ is  dense in $\h$ with respect to $\tau_w$ and $\bar p$.
\end{theorem}
{\it Proof}. (1) The fact that $\overline{\ap}^\tau = \ap$ may be proved in
the same way as in [BDK]. (Recall that by the convention from Section 1 we
take the part of the $\tau$-closure of $\ap$ that lies in $\h$).
Furthermore, since $\tau\succ \tau''$, we have $\overline{\ap}^{\tau''}
\supseteq \overline{\ap}^\tau $. To prove (1), we need to show that the
above inclusion is in fact equality. We will use the equivalence of
$\tau''$ and $\bar\tau$ (see (\ref{bartau})). Let $S\in
\overline{\ap}^{\bar\tau}$ and assume $S$ has a point $x_0$ of period $n$.
Then $(x_0, x_1,...,x_{n-1})$ is a finite $S$-periodic orbit, where $x_i =
S^ix_0$ and $S^nx_0 = x_0$. Take $\mu_i = \de_{x_i},\ i= 0,...,n-1$, and
consider an arbitrary homeomorphism $T$ from $\overline V= \overline V(S;\
\mu_0,...,\mu_{n-1};\ \e )$. It follows that $T$ has the same periodic
orbit $(x_0,...,x_{n-1})$. To see this, assume that $Tx_0 \neq Sx_0$. Then
there exists a clopen set $ F$ containing $Tx_0$ which does not contain
$Sx_0$. By (\ref{bartau}) this contradicts the fact that $T \in \overline
V$. Similarly, one can show that $Tx_i = Sx_i$ for $i =1,...,n-1$. Thus,
every such homeomorphism $T$ has a periodic orbit. This contradicts our
assumption that there exists some $S \in \overline{\ap}^{\tau''} \setminus
\ap$.

We remark that $\ap$ is not closed in $\tau_w$. Indeed, one can easily find
a sequence of aperiodic homeomorphisms that converges  to the identity map
in $\tau_w$.
\smallskip

(2) To prove that $\ap$ is dense in $(\h,\tau_w)$ (and hence in $(\h,\bar
p)$), it suffices to show that each simple homeomorphism can be
approximated by an aperiodic homeomorphism. We use the above notation from
the definition of simple homeomorphisms.

Let $S$ be a simple homeomorphism and let $\e >0$. Denote by
$$
\tilde F_j = \bigcup_{i=0}^{r_j-1} S^iF_j,\ \ j = 1,...,k,
$$
$$
\tilde C_s = \bigcup_{n\in \Z} S^nC_s \cup \{y^+_s,...,S^{q^+_s-1}y^+_s\}
\cup  \{y^-_s,...,S^{q^-_s-1}y^-_s\},\ s=1,...,l.
$$
The sets $\tilde F_1,...,\tilde F_k$ and $\tilde C_1,..., \tilde C_l$ are
clopen, disjoint, $S$-invariant and by (\ref{simple})
$$
\Om = \bigcup_{j=1}^k \tilde F_j \cup \bigcup_{s=1}^l \tilde C_s.
$$

Given $\e$, we will find an aperiodic homeomorphism $T$ such that
\begin{equation}\label{appr}
d_w(S,T)= \sup_{x\in \Om} d(Tx,Sx) +  \sup_{x\in \Om} d(T^{-1}x,S^{-1}x)
<\e.
\end{equation}
To do this, it suffices to find aperiodic homeomorphisms $P_j : \tilde F_j
\to \tilde F_j$ and $R_s : \tilde C_s \to \tilde C_s$ satisfying
(\ref{appr}) on the sets $\tilde F_j$ and $\tilde C_s$ for all $j,s$. To
construct $P_j, j=1,...,k,$ we divide the $S$-tower
$(F_j,...,S^{r_j-1}F_j)$ into finitely many clopen subtowers
$(F_{jm},...,S^{r_j-1}F_{jm}),\ m=1,...,m_j,$ such that ${\rm diam}
(S^iF_{jm}) <\e$ for all $i$ and $m$. Let $P_j(m)$ be an aperiodic
homeomorphism of $F_{jm}$. Define $P_jx = Sx$ for $x\in
\bigcup_{i=1}^{r_j-1}S^iF_{jm}$ and $P_jx = SP_j(m)x$ for $x\in F_{jm}$,
$m=1,...,m_j$. By construction, $P_j$ maps $\tilde F_j$ onto itself and
$d_w(P_j, S) <\e$ on each $\tilde F_j,\ j=1,...,k.$

Fix some $s\in \{1,...,l\}$. To construct an aperiodic homeomorphism $R_s$
of $\tilde C_s$ such that $d_w(S,R_s) < \e$, we will use the following
property (*): given a proper clopen subset $A$ of a Cantor set $Z$, one can
find a sequence of disjoint clopen sets $A_1 = A, A_2,...$ in $Z$ and a
homeomorphism $R: Z\to Z\setminus A$ such that (i) the set $Z\setminus
\bigcup_{j=1}^\infty A_j$ is uncountable, (ii) $RA_j = A_{j+1},\
R(Z\setminus \bigcup_{j=1}^\infty A_j)= X\setminus \bigcup_{j=1}^\infty
A_j$ and $R$ is aperiodic on $Z\setminus \bigcup_{j=1}^\infty A_j$.

Let $a$ be the minimum of distances between the points $\{S^iy^+_s,\
S^jy^-_s : i=0,1,...,q^+_s-1; j=0,1,...,q^-_s-1\}$. Given $0 <\e < a/2$, we
can find $n_0$ such that ${\rm dist}(S^nC_s, S^ny^+_s) < \e/4$ and ${\rm
diam} (S^nC_s) <\e/4$ for $n \geq n_0$. Without loss of generality we can
assume that $n_0 \equiv 0 \ \mod (q^+_s)$.

Denote by
$$
B_p = \bigcup_{i=0}^\infty S^{n_0+p +iq^+_s}(C_s) \cup \{S^py^+_s\},\ \
p=0,1,...,q^+_s -1.
$$
The set $B_p$ is clopen and ${\rm diam} (B_p) < \e/2$ for each $p$. Observe
that $SB_p = B_{p+1}, p=0,...,q^+_s -2$, and $SB_{q^+_s-1} = B_0 \setminus
A$ where $A = S^{n_0}C_s$. Now we can apply property (*) for $Z= B_0$.
Choose an infinite sequence of clopen sets $A_1 = A, A_2, A_3,...$ such
that every $A_i$ is a subset of $B_0$ and the set $B_0 \setminus
\bigcup_{i=1}^\infty A_i$ is uncountable. Take a homeomorphism $R_0$
(defined on $B_0$ only) which maps $B_0$ onto $B_0\setminus A$ and
satisfies the condition:
$$
 R_0(B_0 \setminus \bigcup_{i=1}^\infty A_i) = B_0 \setminus
\bigcup_{i=1}^\infty A_i,\ \ R_0A_i = A_{i+1},\ \ i=1,2,... .
$$
Let now $R_i$ be a homeomorphism defined on $B_{i-1}$ such that $R_iB_{i-1}
= B_i,\ i=1,...,q^+_s-1$. Define the homeomorphism $R^+_s : D_s \to D_s$
where $D_s =\bigcup_{p=0}^{q^+_s-1} B_p$ as follows: $R^+_sx = R_px$ for
$x\in B_p,\ p=0,...,q^+_s-2$ and $R^+_sx = R_0(R_1^{-1}\cdots
R_{q^+_s-1}^{-1})x$ for $x \in B_{q^+_s-1}$. We see that $R^+_sB_{q^+_s-1}
= B_0\setminus A$ and therefore $SB_p = R^+_sB_p$ for all $p$. It follows
that $d(R^+_sx,Sx) + d(R^{+{-1}}_sx,S^{-1}x) < \e$ for $x\in D_s$.

Replacing $S$ by $S^{-1}$ in the above construction, we can similarly
define a homeomorphism $R_s^-$ which acts only on the clopen set $D'_s=
\bigcup_{r=0}^{q^-_s-1} B'_r$ where
$$
B'_r = \bigcup_{i=0}^\infty S^{-m_0- r -iq^-_s}(C_s) \cup \{S^{-r}z^-_s\},
\ \ r= 0,1,...,q^-_s-1,
$$
Here $m_0\equiv 0 \mod (q^-_s)$ is defined analogously to $n_0$ and $z^-_s
= S^{q^-_s -1}y^-_s$. It can be easily checked that $d_w(S,R^-_s) < \e$ on
the set $D'_s$.

Finally, we define $R_s : \tilde C_s \to \tilde C_s,\ s=1,...,l$ as follows
$$
R^-_sx = \left\{
\begin{array}{lll}
    Sx, & \hbox{$x\in \bigcup_{i=-m_0+1}^{n_0-1}S^iC_s$} \\
    R^+_sx, & \hbox{$x\in B_s$} \\
    R^-_sx, & \hbox{$x\in B'_s$}. \\
\end{array}
\right.
$$
Thus, the aperiodic homeomorphism $T$ defined by $P_j$ and $R_s$\
($j=1,...,k; s=1,...,l$) satisfies (\ref{appr}). \hfill$\square$
\\

We note that every simple homeomorphism has a nontrivial periodic part $Z =
\bigcup_{j=1}^k\bigcup_{i=0}^{r_j-1} S^iF_j$. Therefore the two dense
subsets,  $\ap$ and ${\cal S}$, are disjoint in $\h$.
\smallskip

Let $\per_0$ be the subset of $\per$ consisting of all homeomorphisms with
finite period, that is $P\in \per_0$ if and only if there exists $m\in \N$
such that $P^mx =x$ for all $x\in \Om$. This means that $\Om$ can be
decomposed into a finite union of clopen sets $\Om_p$ such that the period
of $P$ at each point from $\Om_p$ is exactly $p$. By $\per_p$, we denote the
subset of $\per_0$ consisting of homeomorphisms with $\Om_p=\Om$. Such
homeomorphisms are called $p$-periodic. Clearly, the set of simple
homeomorphisms, ${\cal S}$ contains $\per_0$.

Let $P\in \per_p$, then any $P$-orbit consists of $p$ different
points. A subset $E\subset \Om$ is called {\it fundamental} for
$P$ if $(E,P(E),...,P^{p-1}(E))$ is a partition of $\Om$.

\begin{lemma}\label{per} Let $\Om$ be a Cantor set  and let $P$ be a
$p$-periodic homeomorphism. Then there exists a clopen $P$-fundamental
subset $E\subset \Om$.
\end{lemma}

\noindent{\it Proof} (suggested by B.~Weiss). Let $d$ be a
metric on $\Om$ compatible with the clopen topology. We note that
there exists some $c>0$ such that $d(x,P^i(x))>c$ for all $x\in \Om$ and
all $i=1,..., p-1$. Indeed, let us fix some $i< p$ and assume that for any
$n\in {\N}$ there exists $x_n\in X$ such that $d(x_n,P^i(x_n))<\frac1{n}$.
Take a convergent subsequence $\{x_{n_k}\} \subset \{x_n\}$ such that
$x_{n_k}\to x_0$ as $k\to \infty$. Then $P^i(x_{n_k})\to P(x_0)$, and
therefore $d(x_0,P^i(x_0))=0$. This contradicts the assumption that $P\in
\per_p$.

Now let $(A_1,A_2,...,A_n)$ be a partition of $\Om$ into clopen sets such
that
\begin{equation}\label{diam1}
{\rm diam}(A_i)\le \frac{c}{2} \qquad \mbox{for all $i=1,...,n$}.
\end{equation}
Define $E_1=A_1$, and for $i=2,3,...,n$, set inductively
\begin{equation}\label{E_i}
E_i=E_{i-1}\cup(A_i \setminus O_{P}(E_{i-1}))
\end{equation}\
where $O_P(F)=\bigcup_{i=0}^{p-1}P^i(F)$.

We first prove that
\begin{equation}\label{disj}
E_k\cap P^i(E_k)=\emptyset,\ \ \ k=1,...,n,\ i=1,...,p-1.
\end{equation}
Clearly, (\ref{disj}) is true for $k=1$. Assume that this relation
is valid for $E_{k-1}$. Then it follows from (\ref{diam1}) and (\ref{E_i})
that
$$
E_k\cap P^i(E_k)=[E_{k-1}\cap P^i(E_{k-1})]\cup [E_{k-1}\cap
(P^i(A_k)\setminus O_P(E_{k-1})]\cup
$$
$$
\cup [P^i(E_{k-1})\cap (A_k \setminus O_P(E_{k-1}))]\cup
[(A_k \setminus O_P(E_{k-1}))\cap (P^i(A_k)
\setminus O_P(E_{k-1}))]=\emptyset.
$$
Next, we show that
\begin{equation}\label{soder}
\bigcup_{i-0}^{p-1}P^i(E_k)\supset \bigcup_{j=1}^k A_j,\ \ k=1,2,...,n.
\end{equation}
Again assume that (\ref{soder}) is proved for $E_{k-1}$. Then
$$
\bigcup_{i=0}^{p-1}P^i(E_k)=\bigcup_{i-0}^{p-1}P^i(E_{k-1})\cup
\bigcup_{i=0}^{p-1}(P^i(A_k)\setminus O_P(E_{k-1}))
$$
$$
=\bigcup_{i=0}^{p-1}P^i(E_{k-1})\cup (A_k \setminus O_P(E_{k-1}))\cup
\bigcup_{i=0}^{p-1}(P^i(A_k) \setminus O_P(E_{k-1})).
$$
The first term contains $A_1\cup\cdots\cup A_{k-1}$ by assumption. The
first and second terms together contain $A_k$.

Thus, it follows from (\ref{disj}) and (\ref{soder}) that for the
clopen set $E=E_n$ the orbit $O_P(E)$ consists of pairwise
disjoint sets and $O_P(E) =\Om$. \hfill$\square$
\\

It follows immediately from Lemma \ref{per} that for every $P\in \per_0$
there exists a $P$-invariant partition $(\Om_1,...,\Om_m)$ of $\Om$ into
clopen subsets such that $\Om_i = \bigcup_{j=0}^{k_i-1}P^{j}E_i$ where
$E_i$ is a fundamental clopen subset for $P$ on $\Om_i$ and $k_i$ is the
period of $P$ on $\Om_i$.
\smallskip

We will now consider the closure of $\per$ in $\h$ with respect to
both $\tau$ and $\tau_w$. Firstly,  we show that $\overline{\per}^{\tau_w}$
is a proper subset in $\h$. This means that there are homeomorphisms which
cannot be approximated by periodic homeomorphisms in $\tau_w$.

We call a homeomorphism $T\in \h$ {\it dissipative} if there exits a clopen
set $F\subset \Om$ such that either $TF\varsubsetneqq F$ or
$F\varsubsetneqq TF$. Clearly, dissipative homeomorphisms exist in
$\h$ since any two clopen sets are homeomorphic.

\begin{proposition}\label{per0}  $(1)$ The set $\overline{\per}^{\tau_w}$
is a proper subset in ($\h,\tau_w)$: In fact, if $T$ is a dissipative
homeomorphism of $\Om$, then $T\notin\overline{\per}^{\tau_w}$.\\
$(2)$ $({\cal S}\setminus \per_0)\cap \overline{\per}^{\tau_w} =
\emptyset$; in other words, if a simple homeomorphism $S$ has an
aperiodic part, then it cannot be approximated in $\tau_w$ by
pointwise periodic homeomorphisms.\\
\end{proposition}

\noindent
{\it Proof}. (1)  Take a dissipative homeomorphism $T$ and let
$F$ be a clopen subset such that $TF \varsubsetneqq F$. We will show that
the neighborhood $W(T; F)$ does not contain any homeomorphism from
$\per$. Assume that this is false and let $P\in \per$ be such that
$PF = TF$.  Then  $P^nF \varsubsetneqq \cdots \varsubsetneqq PF\ (=\
TF)  \varsubsetneqq F$ for any $n$. It follows that there are points
from $F$ with infinite orbits, and this contradicts the pointwise
periodicity of $P$.

(2) It suffices to show that each $S$ from  $({\cal S}\setminus \per_0)$ is
dissipative. We use  notation from the definition of simple homeomorphisms.
Note that it follows from decomposition  (\ref{simple}) that every closed
set $E_s = (\bigcup_{n=0}^\infty S^nC_s)\ \cup\
[(y_s^+,...,S^{q^+_s-1}y_s^+) ]$ is, in fact, clopen because $\Om \setminus
Z$ is a finite disjoint union of closed sets. Clearly, $SE_s\varsubsetneqq
E_s$ and the result follows from (1). \hfill{$\square$}
\\

In Section 4 we will strengthen the above result and give a
complete description of the set $\overline{\per_0}^{\tau_w}$.
\smallskip

Let $T$ be a minimal homeomorphism of a Cantor set $\Om$. We consider the
full group $[T]$ and the {\it topological full group} $[[T]]$ of
homeomorphisms generated by $T$ (see [BK1, GPS2, GW] for details). Recall
that a homeomorphism $\g\in \h$ belongs to $[[T]]$ if and only if $\g x =
T^{m_\g(x)}x$ where $x \mapsto m_\g(x)$ is a continuous function $\Om\to \Z$
.

It was shown in [BK1] how one can use Kakutani-Rokhlin
partitions to describe the structure of homeomorphisms from $[[T]]$. Here
we recall some facts that will be used later on.

Let $(A_n)$ be a sequence of clopen subsets of $\Omega$ such that
$A_n\supset A_{n+1},\ n\in \N,$ and $\bigcap_n A_n$ is a singleton in
$\Omega$. Given $T$ and $A_n$, we can produce a Kakutani-Rokhlin partition
$\xi_n$ of $\Omega$ which is determined by the function of  first return
to $A_n$ under the action of $T$ ([P, HPS]). The partition $\xi_n$ consists
of a finite collection of $T$-towers $\xi_n(v), v\in V_n:$
$$
\xi_n(v)=\{ P_n^i(v):=T^iP_n(v)\ |\ i=0,1,...,h_n(v)-1\}
$$
where $P_n^0(v) = P_n(v)$. We note that $(A_n)$ can be
also chosen such  that $\xi_{n+1}$ refines $\xi_n$ and $\bigcap_n\xi_n$
generates the clopen topology on $\Omega$.  Moreover one can assume that
${\rm diam}(A_n) \to 0$ as $n\to \infty$.

Suppose $\g\in [[T]]$. Then there exists $N\in \N$ such that $E_i=\{x\in
\Om:  \g x=T^ix\}, \ \-N \le i\le N,$ is a clopen finite partition
of $\Om$ (some of $E_i$'s may be empty). Note that for sufficiently large
$n$, each set $E_i$ becomes a $\xi_n$-set, that is $\g P^i(v_n) =
T^lP^i(v_n)$ for some $l = l(i,v_n)$. Moreover, we may suppose that $N<
2^{-1}h_n$, where $h_n=\min(h_n(v) : v\in V_n)$.

We will now commence  the study of the periodic and aperiodic parts of
the topological full group with respect to the uniform and weak topologies
$\tau$ and $\tau_w$.

Let us denote by $\per_0(T)$ the set $\per_0\cap [[T]]$.

\begin{theorem}\label{per1} $\Min \subset \overline{\per_0}^\tau$ and
$\Min \subset \overline{\per_0}^{\tau_w}$. More precisely, let $T$ be
a minimal homeomorphism of $\Om$, then: $(1)$ given a neighborhood
$U(T; \mu_1,...,\mu_m; \e)$, there exists a periodic homeomorphism $P\in
\per_0(T)$ such that $P\in U(T; \mu_1,...,\mu_m; \e)$; $(2)$ given $\e >0$
there exists $Q\in \per_0(T)$ such that $d_w(T, Q) < \e$.
\end{theorem}

\noindent {\it Proof}. We will prove (1); assertion (2) can be proved
similarly.

Every measure $\mu \in {\cal M}_1(\Om)$ has an at most countable set of
points of positive measure; denote it by $\{x_{\mu}(k)\}$. Given
$\mu_1,...,\mu_m$ and $\e>0$, find a finite set $Y = \{ x _{\mu_i}(k) \ :\
i=1,...,m,\ k\in I(\mu_i)\subset \N\}$ where the finite subset $I(\mu_i)$ is
determined by the condition \begin{equation}\label{meas} \sum_{k\notin
I(\mu_i)} \mu_i(\{x_{\mu_i}(k)\}) \le \frac{\e}{3},\ \ i=1,..., m.
\end{equation}
Let $Y = (y_1,...,y_N)$. Choose a point $\bar x\in \Om \setminus Y$ such
that $T\bar x = \bar y$ does not belong to $Y$. By [HPS], we can
find a sequence $(\xi_n), \ \xi_n=\{T^i D_j(n)\ :\ 0\le i\le k(j,n) -1,\
j\in K(n)\},\ |K(n)| < \infty$, of Kakutani-Rokhlin partitions satisfying
the following conditions:

(i) $\xi_{n+1}$ refines $\xi_n$ and for the base $B(n) =\bigcup_jB_j(n)$, one
has $B(n+1) \subset B(n)$;

(ii) $(\xi_n)$ spans the clopen topology on $\Om$;

(iii) $\bigcap_n B(n) =\{\bar y\}$ and $\bigcap_n C(n)= \{\bar x\}$ where
$C(n) = \bigcup_jT^{k(j,n)-1}D_j(n)$.

\noindent
Let $n_0$ be sufficiently large so that $B(n_0)\cap Y = C(n_0)\cap
Y = \emptyset$ and
\begin{equation}\label{<e}
\mu_i(B(n_0)) < \frac{\e}{2},\ \ \ \ \ \mu_i(C(n_0)) < \frac{\e}{2},\ \ \
i=1,...,m.
\end{equation}
The sets $B(n_0)$ and $C(n_0)$ may contain points of positive
measure $\mu_i$ but, by (\ref{meas}), the total contribution of
these points to the measures of either of the sets is less than
$\e/3$.

For every $T$-subtower $\xi^j_{n_0} = \{T^iD_j(n_0)\ : \ 0\le i\le k(j,n_0)
-1\},\ j\in K(n_0),$ we define a periodic homeomorphism $P_j(n_0)$:
\begin{equation}\label{P}
P_j(n_0)x =\left\{ \begin{array}{ll} Tx, & {\rm if}\ x\notin
T^{k(j,n_0)-1}D_j(n_0)\\ T^{-k(j,n_0)+1}x, & {\rm otherwise}
\end{array} \right.
\end{equation}
We define the periodic homeomorphism $P$ as follows: $Px =P_j(n_0)x$
if $ x\in \xi_{j,n}$. We get from (\ref{P}) that $P\in [[T]]$ and
$E(P, T) = B(n_0)\cup C(n_0)$. Thus, by (\ref{<e}), we obtain that
$P\in U(T; \mu_1,...,\mu_m; \e)$.

To prove (2), we observe that ${\rm diam}(B(n))$ and ${\rm diam}(C(n))$
tend to 0 as $n\to \infty$. Therefore, the above method allows us to find
a periodic homeomorphism $Q\in \per_0(T)$ which is $\e$-close to $T$ with
respect to $\tau_w$. \hfill$\square$
\\

Although we have shown that $\per_0$ is not dense in $(\h, \tau_w)$ it is
interesting to decide whether $\per_0(T)$ is dense in $[T]$ with respect
to $\tau$ and $\tau_w$.

\begin{theorem}\label{per2}
Let $(\Om, T)$ be a Cantor minimal system, then:
$(1)$ $\overline{\per_0(T)}^\tau= [T]$ and  $(2)$
$\overline{\per_0(T)}^{\tau_w} \supset  [[T]]$.
\end{theorem}

\noindent {\it Proof}. Case (1) will be considered in detail, case (2)
is similar.

We use here notation from the preceding
proof. Let $\g \in [[T]]$. Then there exists  $K\in \N$ such that for
all $ i\in [-K,K], $ the clopen sets $E_i =\{x\in \Om : \g x= T^i x\}$
constitute a partition $\eta= \eta(\g)$ of $\Om$. We first prove that for
any neighborhood $U_\g = U(\g; \mu_1,...,\mu_m;\e)$ there exists a periodic
homeomorphism $P\in {\cal P}$ such that $P\in U_\g$. We apply the
method used in the proof of Theorem \ref{per1}. Let $Y$ and $(\xi_n)$ be as
above. In addition to (i) - (iii), we may assume that $(\xi_n)$ satisfies
the following conditions (see [BK1] for details):

(iv) the height $k(j,n)$ of every $T$-subtower
$\xi^j_n$ approaches to infinity as $n\to \infty$;

(v) $\xi_n$ refines
$\eta$, i.e. every $E_i$ is a union of atoms from $\xi_n$; in particular,
for every element $D\in \xi_n$, $\g x = T^i x,\ x\in D$, where $i=i(D)$.

Take $M\in \N$ and choose $n_1$ so large that $\min\{k(j,n) : j\in K(n)\}
\ge K(M+2)$ for all $n\ge n_1$. Let
$$
Z =(\bigcup_{i=0}^{K-1}T^iB(n_1))\cup (\bigcup_{i=0}^{K-1}T^{-i} C(n_1)),
$$
where $B(n_1)$ and $ C(n_1) $ are the base and top of $\xi_{n_1}$,
respectively. It follows from (i) - (v) that the $\g$-orbit of any atom $D$
of the partition $\xi_n$ meets $Z$ at least once. Furthermore, by the same
reasoning as in the proof of Theorem \ref{per1}, we can assume that $n_1$ is
chosen sufficiently large that $\mu_i(Z) <\e,\ i=1,...,m$.

Now fix a $T$-subtower $\xi^j_{n_1}$, consisting of sets $T^iD_j(n_1)=
D(i,j)$. Define a periodic homeomorphism $P(j,n_1) : \xi^j_{n_1} \to
\xi^j_{n_1}$ from $[[T]]$ as follows. Let $D(K,j)$ be the first atom (with
respect to the natural order on $\xi^j_{n_1}$) that does not belong to $Z$.
Consider the sets $\g^pD(K,j),\ p=0,...,L$, where $ \g^L(D(K,j))\subset Z $
and $ \g^p D(K,j) \cap Z =\emptyset$. Define $P(j,n_1)x =\g x$ on $
\bigcup_{0\le p< L}\g^p (D(K,j)) $ and $P(j,n_1)x = \g^{-L} x $ on
$\g^L(D(K,j))$. Let $D(i_1,j)$ be the first atom in $\xi_{n_1}^j$ where $
P(j,n_1) $ has not been defined. We extend the definition of $ P(j, n_1) $
on a finite piece of the $\g$-orbit outgoing from $ D(i_1,j)$ that does not
meet $Z$. Repeating this construction we eventually define $P(j,n_1)$ for
all $D(i,j)$ with $K\le i < k(j,n_1) -K +1$ and for some $ D(i,j) $ from
$Z$. We set $P(i,n_1)$ to be the identity map for the remaining part of
atoms of $Z$. Let now $Px = P(j,n_1)x$ if $x$ is in $\xi^j_{n_1}$. Then $P$
is a periodic homeomorphism from $[[T]]$ whose period is at most $M$ at
every point. By construction, $E(P,\g)\subset Z $, that is $P\in U_\g$. To
complete the proof of (1), use the argument of [Theorem 4.5, BK1] where the
density of $[[T]]$ in $[T]$ was established.

For case (2), first observe that ${\rm diam}(Z)$ can be made arbitrarily
small by choosing $n$ sufficiently large. Then use the same method to prove
that each $\g \in [[T]]$ can be approximated in $\tau_w$ by homeomorphisms
from $\per_0(T)$. As was shown in [BK1], the $\tau_w$-closure of $[[T]]$
does not, in general, contain $[T]$. \hfill$\square$
\\

Our next goal is to describe periodic and aperiodic homeomorphisms $\g$ from
$[[T]]$ where $T$ is a minimal homeomorphism of $\Om$. We will refine the
results proved in [BK1], describing all possible types of $\g$-orbits.

Fix some $\xi = \bigcup_{v\in V}\xi(v)$ from the sequence $(\xi_n)$ of
Kakutani-Rokhlin partitions built by $T$ and a refining sequence of clopen
subsets $(A_n)$ (see above). Given $\xi$ and $T$, define two partitions
$\alpha$ and $\alpha'$ of $V$: we say that $J$ is an atom of $\alpha$
if $J$ is the smallest subset of $V$ such that $T(\bigcup_{v\in
J}T^{h(v)-1}D(v))$ is a $\xi$-set. Similarly, $J'\in \alpha'$ if
$T^{-1}(\bigcup_{v\in J'}D(v))$ is a $\xi$-set and $J'$ is the smallest
subset with this property. Obviously, for every $J\in \alpha$ there exists
$J'\in \alpha'$ such that
\begin{equation}\label{J}
T(\bigcup_{v\in J}T^{h(v)-1}D(v)) =\bigcup_{v\in J'}D(v)
\end{equation}
Notice that (\ref{J}) defines a one-to-one correspondence $ m:J\to J'$
between atoms of $\alpha$ and those of $\alpha'$.

For $J\in \alpha$ and $J'\in \alpha'$  define $L_{tj}(J)$ and $L_{bk}(J')$,
$0\le j, k \le \frac{1}{2}h,\ h=\min_{v\in V}h(v)$  as follows:
\begin{equation}\label{lines}
L_{tj}(J)=\bigcup_{v\in J}T^{h(v)-j-1}D(v),\ \
L_{bk}(J')=\bigcup_{v'\in J'}T^kD(v')
\end{equation}
(here $t$ stands for \textquotedblleft top\textquotedblright and $b$ stands
for \textquotedblleft base").  Remark that the indexes $j$ and $k$ in
$L_{tj}(J)$ and $L_{bk}(J')$ indicate the distance of $D^{h(v)-j-1}(v),\
v\in J,$ and $D^k(v),\ v'\in J',$ from the top and from the base of the
corresponding towers.

Since $(\xi_n)$ generates the clopen topology, we note that given $\g \in
[[T]]$ there exists $\xi\in (\xi_n)$ such that for every $ i\in \Z$
the set $\{x\in \Om : \g x = T^i x\},$ is a $\xi$-set. Then for every
$\xi$-atom $D^i(v)$, there exists an integer $l = l(D^i(v)) $ such that
\begin{equation}\label{2.10}
\g x = T^l x,\ \ x\in D^i(v).
\end{equation}
On the other hand, it was proved  in [BK1] that for $\g$ and  $\xi$
as above, the following property holds: if $\g(D^i(v)) = T^l(D^i(v))$ and
$l+i \ge h(v)$ (i.e. $D^i(v)$ goes through the top of $\xi(v)$ under
action of $\g$), then the entire set $L_{tj}(J),\ j = h(v) -i -1$,
containing $D^i(v)$, also goes through the top of $\bigcup_{v\in
J}\xi(v)$. Furthermore, $L_{tj}(J)$ is mapped by $\g$ onto $L_{bk}(J')$
where $J'= m(J)$ and $k$ is uniquely determined by $j,l$. A similar
property holds when $\g(D^k(v')) = T^s(D^k(v')),\ v'\in J'$, and $k + s <
0$. In this case, the set  $L_{bk}(J')$ goes through the base and is
mapped by $\g$ onto some $L_{tj}(J)$.
\smallskip

It turns out that the above result allows us to solve the
inverse problem, that is, to find a finite collection of objects that
uniquely determine a homeomorphism $\g \in [[T]]$. For this, we take the
following data:
\smallskip

(a) a positive integer $N < \frac{1}{2}h$;

(b)  subsets ${\rm I}_L(J),\ {\rm II}_L(J'),\ {\rm I}_A(J),\ {\rm
II}_A(J')$  of $\{ 0,1,...,N\}$ such that  $|{\rm I}_L(J)|=
|{\rm II}_A(J')|,\  |{\rm II}_L(J')|=|{\rm I}_A(J)|$,  $J\in \alpha,\
J'\in \alpha'$;

(c)  one-to-one maps $\rho(J): {\rm I}_L(J)\rightarrow {\rm II}_A(J'),\
\sigma(J'): {\rm II}_L(J')\rightarrow {\rm I}_A(J)$
where $J\in \alpha,\ J'\in \alpha'$, and $m(J)=J'$;

(d)  a one-to-one map  $\pi(v): {\rm I}_{NL}(v)\rightarrow {\rm
I}_{NA}(v) $ where $v\in J\cap J', \ J\in \alpha,\ J'\in \alpha'$, and
the set ${\rm I}_{NL}(v)$ (resp. ${\rm I}_{NA}(v)$) consists of those
$k\in \{ 0,1,...,h(v)-1\}$ such that $k\notin  {\rm  II}_L(J')$
(resp. $k\notin {\rm II}_A(J')$) and $h(v) - k -1\notin {\rm I}_L(J)$
(resp. $h(v) - k -1\notin {\rm I}_A(J)$).
\smallskip

In the above notation, the indexes $A,L,NA,NL$ mean the first letters in
words \textquotedblleft arriving", \textquotedblleft leaving",
\textquotedblleft not arriving", \textquotedblleft not leaving".

Using these data, we can define a homeomorphism $\g$
from $[[T]]$ by the following rule:

For every $\xi$-atom  $D^j(v) = T^j(D(v)),\ v\in J\cap J'$,  we have that:
either (1) $j$ is  in ${\rm I}_{NL}(v)$, or (2)  $h(v) - j- 1\in {\rm
I}_L(J)$, or (3) $j\in {\rm II}_L(J')$. (Note that in view of (a) - (d)
only one of the above possibilities can occur). If (3) holds, then
$D^j(v)$ belongs to $L_{bj}(J')$; if (2) holds, then $D^j(v)$ belongs to
$L_{tk}(J)$ where $k = h(v) -j-1$.

According to cases (1) - (3),  we define
\begin{equation}\label{gamma} \gamma(T^jD(v))=T^{\pi (v)(j)}D(v),\ \ \
{\rm if}\ \  j\in {\rm I}_{NL}(v),
\end{equation}
\begin{equation}\label{gamma1}
\gamma(L_{tk}(J))=L_{b,\rho (J)(k)}(J'),\ \ \ {\rm if}\ \   j\in
{\rm I}_L(J),
\end{equation}
\begin{equation}\label{gamma2}
\gamma(L_{bj}(J'))=L_{t,\sigma (J')(j)}(J), \ \ \ \ {\rm if}\ \
j\in {\rm II}_L(J').
\end{equation}
We observe that the image of $D^j(v)$ under the $\g$-action is a
$\xi$-set if $\g$ is defined by (\ref{gamma}) and is no longer a $\xi$-set
if $\g$ is defined by (\ref{gamma1}) and  (\ref{gamma2}).

Denote by $\G(\xi)$ the set of all homeomorphisms that can be constructed
from the data  (a) - (d) by (\ref{gamma}) -- (\ref{gamma2}).
It was proved in [BK1] that $\G(\xi_n)$ is an increasing sequence of
subsets in $[[T]]$ and, for every $\g\in [[T]]$, there exists $\xi$ from
$(\xi_n)$ such that $\g\in \G(\xi)$.
\smallskip

Fix a homeomorphism $\g\in [[T]]$.  Let $\xi \in (\xi_n)$ satisfy
(\ref{2.10}).  Then $\g$ determines the subsets   ${\rm I}_L(J),\ {\rm
II}_L(J'),\ {\rm I}_A(J),\ {\rm II}_A(J')$  and maps   $\rho(J),
\sigma(J')$ such that  (\ref{gamma}) -- (\ref{gamma2}) hold.  Let us
consider $\g$-orbits in terms of these subsets and maps.

First suppose that $\g\in [[T]]$ is periodic. There are two possible types
of periodic behaviour for $\g$.

 $1^{\rm st}$ case. We start with some $D^{j_0}(v)$ where $j_0\in {\rm
I}_{NL}(v)\cap {\rm I}_{NA}(v)$. Suppose that the $\g$-orbit of
$D^{j_0}(v)$ does not leave $\xi(v)$. By (\ref{gamma}), this
means that the entire sequence $(j_k)_{k=0}^s, \ j_k=\pi(v)(j_{k-1})$,
belongs to ${\rm I}_{NL}(v)\cap {\rm I}_{NA}(v)$ and $j_s=j_0,\ j_k\neq
j_0,\ k < s$. Let
\begin{equation}\label{eta1}
\eta_1 =\{ T^{j_k}D(v)\ |\
k=0,1,...,s-1\}=\{\gamma^i(T^{j_0}D(v))\ |\ i=0,1,...,s-1\}
\end{equation}
be the  $\g$-orbit where $\gamma(T^{j_k}D(v))=T^{j_{k+1}}(D(v))$.
Then, it follows that the $\g$-orbit of $D^{j_0}(v)$ returns
to this set and because the orbit is a part of $T$-tower, we get that $\g$
is periodic.

$2^{\rm nd}$ case.  The other type of cyclic $\g$-orbit has the following
structure. Fix some  $J\in \alpha$ and let $J'=m(J)$. Suppose that
that $j_0\in {\rm I}_L(J)\cap  {\rm I}_A(J), \ j_1=\rho(J)(j_0)\in {\rm
I}{\rm I}_A(J')\cap {\rm II}_L(J'), \ j_2=\sigma(J')(j_1) \in
{\rm I}_A(J) \cap {\rm I}_L(J),...,j_{2s-1} = \rho(J)(j_{2s-2})\in {\rm
I}{\rm I}_A(J')\cap {\rm II}_L(J')$ and $\sigma(J')(j_{2s-1})=j_0$.  By
 (\ref{gamma1}) and (\ref{gamma2}), this case corresponds to the following
periodic $\g$-orbit
\begin{equation}\label{eta2}
\eta_2=\{\gamma^i(L_{tj_0}(J))\ |\ i=0,1,...,2s-1\}
\end{equation}
where $\g(L_{tj_0}(J)) =  L_{bj_1}(J'),\  \g(L_{bj_1}(J')) = L_{tj_2}(J))$
etc. until  it returns to $L_{tj_0}(J)$.

It follows from the above proof that $\g$ belongs to $\per_0$.
We summarize the above observations in the following statement.

\begin{lemma}\label{percase} Let $\g$ be a periodic homeomorphism from
$[[T]]$. Then there exists $\xi$ from $(\xi_n)$ such that every periodic
$\g$-orbit either has the form $(\ref{eta1})$ or the form $(\ref{eta2})$.
\end{lemma}

We can study $\g$-orbits for an aperiodic $\g \in [[T]]$ in a similar
manner. Given $\e >0$ and $\g\in [[T]]_{ap}$  where $[[T]]_{ap} =
[[T]]\cap\ap$,  we take $\xi$ from the sequence $(\xi_n)$
such that $\g\in \G(\xi)$.   Since ${\rm diam}(A_n)\to 0$ as $n\to
\infty$, we can assume that
\begin{equation}\label{diam}
{\rm diam}\left(\bigcup_{J\in \alpha}\bigcup_{j=0}^N L_{tj}(J)\right) +
{\rm diam}\left( \bigcup_{J'\in \alpha'}\bigcup_{j=0}^N L_{bj}(J') \right)
< \e.
\end{equation}
Considering  $\g$-orbits of atoms of $\xi$, one can find a finite
partition $\zeta$ of  $\Om$ into clopen subsets, $\g$-towers, such that
those towers have their bases and tops into the sets  $\bigcup_{j=0}^N
L_{bj}(J')$ and  $\bigcup_{j=0}^N L_{tj}(J)$. Then we can construct an
odometer which is, by (\ref{diam}), $\e$-close to $\g$ in $\tau_w$ (see
e.g. the proof of Theorem 4.3). The case of the topology $\tau$ is
considered similarly. We should note that given $\e >0$ and
$\mu_1,...,\mu_n \in {\cal M}_1(\Om)$ we can chose $\xi$ such that  for all
$i$,
$$
\mu_i\left(\bigcup_{J\in \alpha}\bigcup_{j=0}^N L_{tj}(J)\right) +
\mu_i\left( \bigcup_{J'\in \alpha'}\bigcup_{j=0}^N L_{bj}(J') \right)
< \e.
$$
This relation guarantees that the odometer which we have found is also
$\e$-close to $\g$ in $\tau$.

Using this ideas, it is now easy to prove the following theorem. We leave
details to the reader. Another proof of the first relation of the theorem
can also be obtained from Theorem \ref{mov} (see Remark \ref{rem} below).

\begin{theorem}\label{ap-gamma} If $T$ is a minimal homeomorphism of
$\Omega$, then $[[T]]_{ap} \subset \overline{\Min}^\tau$ and   $[[T]]_{ap}
\subset \overline{\Min}^{\tau_w}$ where $[[T]]_{ap} = [[T]]\cap\ap$.
\end{theorem}

Let $T$ be a homeomorphism of a Cantor set $\Om$. Then $T$ can be also
considered as a Borel automorphism of $(\Om, \B)$. Therefore, one can define
two full groups $[T]_C$ and $[T]_B$ where
$$
[T]_C = \{S\in \h : Sx \in \{T^nx\ |\ n\in \Z\} \ \forall x\in \Om\},
$$
$$
[T]_B = \{S\in Aut(\Om,\B): Sx \in \{T^nx\ |\ n\in \Z\} \ \forall x\in
\Om\}.
$$
Here the subindeces $C$ and $B$ correspond to the cases of Cantor and Borel
dynamics. Obviously, $[T]_C \subset [T]_B$ and $[T]_B$ is closed in
$Aut(\Om,\B)$ with respect to $\tau$ [BDK].

If $S\in [T]_C$ (or $S\in [T]_B$), then $S$ generates two partitions
$\pi(S)= (X_n : n\in \Z)$ and $\pi'(S) = (X'_n : n\in \Z)$ of $\Om$ into
closed (Borel) subsets $X_n = \{x\in \Om : Sx = T^nx\}$ and $X'_n = S(X_n) =
T^n(X_n)$. Those homeomorphisms from $[T]_C$, for which the sets $X_n$ are
clopen, form the so called topological full group $[[T]]_C$. It was shown in
Section 1 that $\h$ is non-closed and dense in $Aut(\Om, \B)$. On the other
hand, we proved in [BK1] that for a minimal homeomorphism $T$,
\begin{equation}\label{old}
\overline{[[T]]_C}^\tau \cap \h = [T]_C.
\end{equation}
But the problem of finding the entire closure of $[[T]]_C$ in $Aut(\Om, \B)$
with respect to $\tau$ remained open. We answer this question in the
following theorem.

\begin{theorem}\label{T2.9} Let $T$ be a minimal homeomorphism of a Cantor
set $\Om$. Then $\overline{[[T]]_C}^\tau = [T]_B$.
\end{theorem}

\noindent {\it Proof}. It is clear that $\overline{[[T]]_C}^\tau
\subset[T]_B$. Take a Borel automorphism $S\in [T]_B$. Let $U(S) =
U(S;\mu_1,...,\mu_m; \e)$ be a $\tau$-neighborhood of $S$. To prove the
theorem we need to show that $U(S)$ contains a homeomorphism $R$ from
$[[T]]_C$. By (\ref{old}), it suffices to prove that there exists some $R_1
\in [T]_C\cap U(S)$.

Consider the partitions $\pi(S)$ and $\pi'(S)$ defined above. Choose $n_0$
such that
\begin{equation}\label{finite}
\sum_{|n| >n_0}\mu_i(X_n) < \e/4,\qquad \sum_{|n| >n_0}\mu_i(X'_n) < \e/4, \
\ \ i =1,...,m.
\end{equation}
For each $X_n,\ |n| \le n_0,$ find a closed $F_n \subset X_n$ such that for
all $i=1,...,m$,
\begin{equation}\label{F's}
\mu_i(X_n \setminus F_n) < \frac{\e}{4(2n_0+1)},\ \ \ \mu_i(X_n \setminus
S(F_n)) < \frac{\e}{4(2n_0+1)}.
\end{equation}
Let $F'_n = S(F_n)$. Clearly, $F'_n$ is also a closed subset in $X'_n$. The
collections $(F_n)$ and $(F'_n)$ consist of pairwise disjoint sets. Then
there exist clopen sets $G_n \supset F_n$ and $G'_n \supset F'_n$ such that
$(G_n : |n| \le n_0)$ and $(G'_n : |n| \le n_0)$ are pairwise disjoint
collections of sets.

Let $\xi_k = \bigcup_{v\in V_k}\xi_k(v)$ be a sequence of Kakutani-Rokhlin
partitions built by $T$ and a decreasing sequence of clopen sets (see
above). For every $k$, find partitions $\alpha= \alpha_k$ and $\alpha'=
\alpha'_k$ with atoms $J$ and $J'$ satisfying (\ref{J}). Define a new
partition $\eta_k$ of $\Om$ which consists of the sets $L_{tj}(J),
L_{bj}(J'),\ J\in \alpha, J'\in \alpha', j=1,...,n_0$ (see (\ref{lines}))
and the remaining atoms of $\xi_k$. Clearly, every atom of $\eta_k$ is a
$\xi_k$-set and the sequence $(\eta_k)$ generates the clopen topology on
$\Om$.

Choose $k$ sufficiently large such that
$$
\min_{v\in V_k} h(k,v) > 3n_0,
$$
where $h(k,v)$ is the height of $\xi_k(v)$ and every set $G_n,\ G'_n,\ |n|
\le n_0$, is an $\eta_k$-set. Without loss of generality, we can assume that
if $E$ is an atom of $\eta= \eta_k$ such that $E\subset G_n$ (or $E\subset
G'_n$), then $E\cap F_n \neq \emptyset$ (or $E\cap F'_n \neq \emptyset$),
$|n|\le n_0$.

Fix some $G_n$. By construction, $G_n$ is a union of atoms $D_1,...,D_p$ of
$\xi_k$ and sets $L_{tj}(J), L_{bj}(J')$ and each of these sets intersects
$F_n$. Define $\overline{S}x = T^nx$ if $x\in \bigcup_{s=1}^pD_s$. We also
define $\overline{S}x = T^nx$ if $x \in L_{tj}(J)$ (or $x\in L_{bj}(J')$)
and the set $L_{tj}(J)$ ($L_{bj}(J')$) goes through the top (base) of
$\xi_k$ under the action of $T^n$. If $L_{tj}(J)$ (or $L_{bj}(J')$) is a
subset of $G_n$ which does not go through the top (base) under $T^n$, we
discard from $L_{tj}(J)$ ($L_{bj}(J')$) those atoms of $\xi_k$ which do not
meet $F_n$ and set $\overline{S}x = T^nx$ for $x$ from the remaining atoms.
In such a way, we have found a clopen subset $\overline{G}_n\subset G_n$ and
a map $\overline{S}$ defined on $\overline{G}_n$. Clearly, $\overline{G}_n
\supset F_n$. Similarly, we define a clopen set $\overline{G'}_n$ such that
$F'_n \subset \overline{G'}_n \subset G'_n$. It can be easily seen that
$\overline{S}\overline{G}_n = \overline{G'}_n$. Indeed, if $E\subset
\overline{G}_n$ and $E\in \eta$, then there exists $y\in E\cap F_n$. Hence
$T^ny \in F'_n$. At the same time, $T^n(E)$ is an $\eta$-set. Thus,
$T^n(E)\subset \overline{G'}_n$. In such a way, we have found a partially
defined homeomorphism $\overline{S}$ such that $\overline{S}(A) =B$ where
$$
 A =\bigcup_{|n|\le n_0} \overline{G}_n\ \ {\mbox and}\ \
B = \bigcup_{|n|\le n_0}\overline{G'}_n.
$$
Since $\overline{S} \in [[T]]_C$ on the clopen set $\bigcup_{|n|\le n_0}
\overline{G}_n$, we get that for any $T$-invariant measure $\nu$
$$
\nu(X \setminus \bigcup_{|n|\le n_0} \overline{G}_n) = \nu(X \setminus
\bigcup_{|n|\le n_0} \overline{G'}_n).
$$
By [GW, Proposition 2.6], there exists a homeomorphism $R'\in [T]_C$ which
maps $X \setminus \bigcup_{|n|\le n_0} \overline{G}_n$ onto $X \setminus
\bigcup_{|n|\le n_0} \overline{G'}_n$. Now define $Rx = \overline{S}x$ if
$x\in A$ and $Rx = R'x$ if $x\in X\setminus A$. Then $R\in [T]_C$ and it
remains to show that $R\in U(S)$. Note that if $x\in F_n$, then $Sx =
\overline{S}x = Rx$. Therefore
$$
E(S,R) \subset (X \setminus\bigcup_{|n|\le n_0}F_n) \cup (X \setminus
\bigcup_{|n|\le n_0}F'_n).
$$
Then, for given measures $\mu_i,\ i=1,...,m$, we have by (\ref{finite}) and
(\ref{F's})
$$
\mu_i(E(S,R)) < \mu_i(X \setminus \bigcup_{|n|\le n_0}F_n) + \mu_i(X
\setminus \bigcup_{|n|\le n_0}F'_n)
$$
$$
= \mu_i(X \setminus\bigcup_{|n|\le n_0}X_n) + \mu_i(X \setminus
\bigcup_{|n|\le n_0}X'_n) + \mu_i(\bigcup_{|n|\le n_0}(X_n \setminus F_n)) +
\mu_i(\bigcup_{|n|\le n_0}(X'_n \setminus F'_n))
$$
$$
< \e/4 +\e/4 + 2(2n_0 +1)\frac{\e}{4(2n_0+1)} = \e.
$$
\hfill$\square$
\\

\begin{remark} \label{sat}
{\rm We note that the $\tau_w$-closures of the full group $[T]$ and the
topological full group $[[T]],\ T\in \ap,$ can be easily found. Indeed, it
was noted $\overline{[T]}^{\tau_w} = \{S\in \h : \mu \circ S =\mu\ \forall
\mu \in M_1(T)\}$ where $M_1(T)$ is the set of $T$-invariant Borel
probability measures [GPS2]. It is not hard to show that $R\in
\overline{[[T]]}^{\tau_w}$ if and only if for every clopen set $E$ there
exists $\gamma \in [[T]]$ such that $RE = \gamma E$. In other words, $R\in
\overline{[[T]]}^{\tau_w}$ if and only if the clopen sets $E$ and $RE$ are
$[[T]]$-equivalent for every $E$. Recall that in [BK1] we defined the
notion of saturated homeomorphisms: a minimal homeomorphism $T$ is called
saturated if any two clopen sets $A$ and $B$ such that $\mu(A) = \mu(B)$
for all $\mu\in M_1(T)$ are $[[T]]$-equivalent. Obviously, every odometer
is a saturated homeomorphism. It was proved there that $T$ is saturated if
and only if $\overline{[[T]]}^{\tau_w} = \overline{[T]}^{\tau_w}$. Now this
result follows easily from the above description of the $\tau_w$-closures}.
\end{remark}

%%%%%%%%%%%%%%SECTION 3%%%%%%%%%%%%%%%%%%%%%%%%%%

\sect{Rank of a homeomorphism}

%%%%%%%%%%%%%%%%%%%%%%%%%%%%%%%%%%%%%%%%%%%%%%%%%%%%%

\setcounter{equation}{0}

The concept of the rank of an automorphism of a standard measure space is an
important invariant in ergodic theory. This notion has been studied in many
papers. M.~Nadkarni [N] has recently defined rank for Borel automorphisms.
Here we consider {\it rank for homeomorphisms} of a Cantor set.

We first recall the definition of {\it odometer} (or {\it adding machine}).
Let $\{\la_t\}_{t=0}^\infty$ be a sequence of integers such that $\la_t\geq
2$. Denote by $p_{-1}=1,\ p_t=\la_0\la_1\cdots \la_t,\ t=0,1,...\ .$ Let
$\Delta$ be the group of all $p_t$-adic numbers; then any element of
$\Delta$ can be represented as an infinite formal series:
$$
\Delta=\{x=\sum_{i=0}^\infty x_ip_{i-1}\ \vert\
x_i\in (0,1,...,\la_i-1)\}.
$$

It is well known that $\Delta$ is a compact metric abelian group. An
odometer, $S$, is the transformation acting on $\Delta$ as follows:
$S x=x+1,\ x\in \Delta$, where $1=1p_{-1}+0p_0+0p_1+\cdots \in \Delta$.
From the topological point of view, $(\Delta, S)$ is a strictly ergodic
Cantor system.

Let
$$
D_0^t=\{x=\sum_{i=0}^\infty x_ip_{i-1}\ \vert\
x_0=x_1=...=x_t=0\}.
$$
We see that the sets $(D_0^t,D_1^t,...,D_{p_t-1}^t),\
D_i^t=S^i(D_0^t), $ form a partition $\xi_t$ of $\Delta$ into clopen
sets. Clearly, $(\xi_t), \ t\ge 0,$ is a refining sequence of
of $S$-towers.  Moreover, $S(\xi)_t = \xi_t$ for every $t$.

We will denote by ${\cal O}d ={\cal O}d(\Om) $ the set of homeomorphisms of
$\Om$ homeomorphic to an odometer $(\De, S)$. Elements from  ${\cal O}d$
will be also called odometers.

\begin{lemma}\label{odom}  Let $T\in \h$ and let $(F_1,...,F_n)$ be a
partition of $\Om$ into clopen sets such that $TF_i = F_{i+1},\ 1\le i\le
n-1$, $TF_n = F_1$. Then there exists an odometer $S$ such that $SF_i =
TF_i$ for all $i$, i.e. $S\in W(T; F_1,...,F_n)$.
\end{lemma}

The proof of the lemma is based on the definition of odometer and left to
the reader.

\begin{definition}\label{rank}
Let $T \in \h$ and let $\F_n,\ n\in \N,$ be a union of $r$ disjoint
$T$-towers consisting of clopen sets, that is
$$
\F_n = \bigcup_{j=1}^r\bigcup_{i=0}^{h_n(i)-1} T^iF_n(j)
$$
where $F_n(j)\in CO(\Om)$ is the base of $j$-th tower and $h_n(j)$
is its height. We say that $T$ has rank at most $r$ if $\F_{n+1}$
refines $\F_n$ and all $\F_n$'s generate the clopen topology on
$\Om$. We say that $T$ has rank $r$ if $T$ has rank at most $r$
but does not have rank at most $r-1$. $T$ has infinite rank if it
does not have rank $r$ for any finite $r$.
\end{definition}

Obviously, every odometer is a homeomorphism of rank one.

We denote the set of all homeomorphisms having rank at most $r$
by $\R(\le r)$ and the set of homeomorphisms of rank $r$  by
$\R(r)$.

\begin{proposition}\label{allspace}
Let $T\in \h$  and suppose $rank(T) = r <\infty$. If $(\F_n)$ is
a sequence of clopen subsets of $\Om$ as in Definition $2.1$, then $\F_n
=\Om$ for all sufficiently large $n$.
\end{proposition}

\noindent {\it Proof}. We note that $\F_n \subset \F_{n+1}$ for all $n\in \N
$ by Definition 2.1. If we assume that $\F_n \neq \Om$ for all $n$, then we
get that $\bigcap_n(\Om\setminus \F_n)\neq \emptyset$. This
contradicts  the assumption that the $\F_n$'s generate the topology
on $\Om$. \hfill$\square$.
\\

Suppose that the clopen set $\F$ is a disjoint union of $r$ disjoint
families of clopen sets: $\F =\bigcup_{i,j}D_{ij}$ where $i= 1,...,r$ and
$j=0,1,...,h(i) -1$. Let $\sigma$ be a total order on the set
$(0,1,...,h(i)-1)$ for $i=1,...,r$. Denote by $\sigma(j)$ the successor of
$j\in (0,...,h(i)-1)$. Define
\begin{equation}\label{zet}
{\cal Z}_{r,\sigma}(\F) = \{R\in \h\ |\ RD_{ij} = D_{i\sigma(j)},\
i=1,...,r,\ j=0,...,h(i)-2 \}.
\end{equation}
Let $S\in {\cal Z}_{r,\sigma}(\F)$. Then $S$ transforms the sets
$(D_{ij} : j=0,...,h(i) -1)$ into an $S$-tower, provided we introduce a
new enumeration such that $SD_{ij}= D_{ij+1}$. Note that
\begin{equation}\label{neighb}
{\cal Z}_{r,\sigma}(\F) = \bigcap_{i=1}^r\bigcap_{j=0}^{h(i)-2} W(S; D_{ij})
= W(S;D_{10},...,D_{1h(1)-2},D_{20},...,D_{rh(r)-2}).
\end{equation}

\begin{theorem}\label{G_d}
For every finite $r$, the set $\R(\le r)$ is a $G_\de$-set
in $\tau_w$-topology and $\R(\le r)$ is an $F_{\sigma\de}$-set in the
topologies $\tau ,\tau'',\bar p$. In particular, $\R(1)$, the set of
homeomorphisms of rank $1$, is a $G_\de$-set in $\tau_w$ and $\R(1)$ is a
$F_{\sigma\de}$ in $\tau, \tau'', \bar p$.
\end{theorem}

\noindent {\it Proof}. Let $(Q_n)$ be a refining sequence of finite
partitions into clopen sets generating the clopen topology. Then for a
finite partition $\F$ into clopen subsets there exists $Q_n$ such that $Q_n$
refines $F$,\ $Q_n\succ F$. One can easily check that
\begin{equation}\label{rank<r}
\R(\le r) = \bigcap_{n=1}^\infty\bigcup_{k=n}^\infty\bigcup_{Q_n\prec
\F\prec Q_k} \bigcup_\sigma {\cal Z}_{r,\sigma}(\F)
\end{equation}
By (\ref{neighb}) and Proposition \ref{property}, every set ${\cal
Z}_{r,\sigma}(\F)$ is clopen in $\tau_w$ and closed in $\tau,\tau'',\bar
p$. The theorem follows. \hfill$\square$

\begin{theorem}\label{r1-min} If $T$ has rank 1, then $T$ is minimal.
\end{theorem}

\noindent {\it Proof}. Let $(\F)_n$ be a refining sequence of $T$-towers
generating the topology on $\Om$. Without loss of generality, we may assume
that every $\F_n$ is a clopen partition of $\Om$. Take some $x\in \Om$ and
let $E$ be a nonempty clopen subset. Find $n$ such that $E$ is an
$\F_n$-set. Then there exists an integer $m$ such that $T^mx\in E$. This
means that $T$-orbit of $x$ is dense in $\Om$.\hfill$\square$

\begin{theorem}\label{r1=od}  $T\in \R(1)$ if and only if $T$ is
topologically conjugate to an odometer. In other words, $\R(1) = \od$.
\end{theorem}

\noindent{\it Proof}. Let $(\F)_n$ be a sequence of $T$-towers corresponding
to a rank 1 homeomorphism $T$. Then
$$
\F_n = \bigcup_{j=0}^{h_n-1} D^n_j, \qquad TD^n_j = D_{j+1}^n, \ 0\le j\le
h_n-2.
$$
Since $\F_{n+1}\succ \F_n$ and every $\F_n$ is a clopen partition of $\Om$,
we see that $h_n$ divides $h_{n+1}$, say $h_{n+1} =\la_nh_n$. We have that
$$
D^n_j = \bigcup_{k=0}^{\la_n-1} D^{n+1}_{kh_n+j},\quad j=0,...,h_n-1.
$$
This proves that $T$ is completely defined, up to conjugacy, by a sequence
of positive integers $(\la_n)$ and therefore $T$ is conjugate to an odometer
by Lemma \ref{odom}. \hfill$\square$
\\

Let $C(T)$ denote the centralizer of a homeomorphism $T$, that is $C(T) =
\{R\in \h : RT = TR\}$. Denote by $Wcl(T)$ the $\tau_w$-closure of the
cyclic group $\{ T^n : n\in \Z\}$. We prove now the weak closure theorem
for homeomorphisms of rank 1.

\begin{theorem}\label{wcl} Let $T\in \h$ be a homeomorphism of rank 1. Then
$C(T) = Wcl(T)$.
\end{theorem}

\noindent {\it Proof}. To show that $C(T) \supset Wcl(T)$, it suffices to
check that $C(T)$ is closed in $\tau_w \sim p$. To do this, we take a
sequence $(S_n) \subset C(T)$ converging to $S$. Since $\h$ is a topological
group in $\tau_w$, then $(S_n, T)\mapsto S_nT$ and $(T,S_n)\mapsto TS_n$ are
continuous and therefore, taking the limit, we get that $ST= TS$.

To see that $Wcl(T) \supset C(T)$, we use the fact that $T$ can be taken as
an odometer by Theorem \ref{r1=od}. Let $R\in C(T)$ and let $\xi=
(F,TF,...,T^{n-1}F)$ be a $T$-tower covering $\Om$ (then $T^nF =F$). Given
$\de> 0$, we can find sufficiently large $n$ such that ${\rm diam}(F)< \de$.
Denote by $\e_0= \min({\rm dist}(F,T^jF) : 0< j <n)$. Then given $0<\e<
\e_0$ find some $\de$ such that $d(Rx,Ry)< \e$ whenever $d(x,y) < \de$ where
$d$ is a metric on $\Om$. It follows that $RF$ is a subset of some $T^iF,
0\le i_0< n$. If we assume that $RF$ is a proper subset of $T^{i_0}F$, then
we come to a contradiction since, in this case, $R(T^jF)$ must be also a
proper subset of $T^{i_0+j}F, 1\le j<n$, where $i_0+j$ is understood by
${\rm mod}\ 0$. In other words, we have shown that $R \in W(T^{i_0};
F,TF,...,T^{n-1}F)$. \hfill$\square$
\\

We note that J.~King proved in [K] the weak closure theorem in the context
of measurable dynamics.

\begin{remark} {\rm  Let $T$ be an odometer and $\xi= (F,TF,...,T^{n-1}F)$
be a $T$-tower covering $\Om$. Then for any $S\in W(T;F,...,T^{n-1}F)$ we
have that $S^i \in W(T^i;F,...,T^{n-1}F),\ i\in \Z$. More general,}
$$
W(T;F,...,T^{n-1}F)^i \subset W(T^i;F,...,T^{n-1}F),\ \ i\in \Z.
$$
\end{remark}

\sect{Minimal and mixing homeomorphisms}

\setcounter{equation}{0}

Let $\Min$ denote the set of minimal homeomorphisms and let $\mix$ be the
set of all mixing homeomorphisms of $\Om$. The following statement shows
that minimality and mixing are not typical properties.

\begin{proposition}\label{minmix}
The following properties hold:\\
$(1)\ \overline{\Min}^\tau \subset \ap,\ \overline{\mix\cap \ap}^\tau
\subset \ap $.\\
$(2)$ For any neighborhood $W =
W({\mathbb I}; F_1,...,F_n),\ \overline{\Min}^{\tau_w} \cap W =\emptyset$
and $\overline{\mix}^{\tau_w} \cap W =\emptyset$.\\
$(3)$ If $R\in {\cal S}$ is a  simple homeomorphism  (in particular, $R$
can belong to $\per_0$), then there exists a neighborhood $W=W(R;
F_1,...,F_n)$ such that $\overline{\Min}^{\tau_w} \cap W =\emptyset,\
\overline{\mix}^{\tau_w} \cap W =\emptyset$.
\end{proposition}

\noindent{\it Proof}. (1) The result follows from Theorem \ref{ap} (1).

(2) Obviously, $\Min$ and $\mix$ do not meet any neighborhood $W =
W({\mathbb I}; F_1,...,F_n)$ of the identity. Since $W$ is a clopen subset
in $(\h, \tau_w)$, we get that $\overline{\Min}^{\tau_w} \subset W^c$ and
$\overline{\mix}^{\tau_w}\subset W^c$. In particular, we see that
$\overline{\Min}^{\tau_w}$ and $\overline{\mix}^{\tau_w}$ are proper subsets
of $\h$.

(3) Since $R \in {\cal S}$, we can find a clopen subset  $E$ of $\Om$
such that the sets $R^iE,\ i=0,1,...,n-1,$ are disjoint and $R^nE =E$.
Denote $\Om_n = \bigcup_{i=0}^{n-1} R^iE$. Then if $S\in W(R;\
E,RE,...,R^{n-1}E)$, then $S(R^iE) = R^{i+1}E$ and therefore $S$ cannot be
mixing. Now let $F$ be a non-empty clopen subset of $E$. The subsets
$F,RF,...,R^{n-1}F$ are still disjoint and if $S\in W(R;
F,RF,...,R^{n-1}F)$, then  we have that $S(\bigcup_{i=0}^{n-1} R^iF) =
\bigcup_{i=0}^{n-1} R^iF$, that is $S$ is not minimal. \hfill$\square$ \\

Since the set ${\cal S}$ of simple homeomorphisms is dense in $(\h,
\tau_w)$, we easily obtain the following result.

\begin{corollary}\label{nowhere}  The
sets $\overline{\Min}^{\tau_w}$ and $\overline{\mix}^{\tau_w}$ are
nowhere dense in $(\h, \tau_w)$. \end{corollary}

\begin{theorem}\label{min} $(1)\ \overline{\od}^\tau =\overline{\R(1)}^\tau
=\overline{\Min}^\tau, \ \ (2)\ \overline{\od}^{\tau_w}
=\overline{\R(1)}^{\tau_w} =\overline{\Min}^{\tau_w} $.
\end{theorem}

\noindent{\it Proof}. It follows from Theorem \ref{r1=od} that we need to
prove the relations $\Min \subset \overline{\R(1)}^\tau$ and $\Min \subset
\overline{\R(1)}^{\tau_w}$ only.  First consider the closure
$\overline{\R(1)}^{\tau_w}$. Let $T$ be a minimal homeomorphism. Take a
sequence of clopen subsets $(A_n)$ such that $A_n \supset A_{n+1}$ and
$\bigcap_{n\ge 1} A_n$ is a singleton. Let $\xi_n$ be the Kakutani-Rokhlin
partition defined as in Section 2 by $T$ with fixed base $A_n$. It is known
that $(A_n)$ may be chosen in such a way that the refining sequence
$(\xi_n)$ generates the topology on $\Om$. Every $\xi_n$ is a finite union
of $T$-towers $\xi_n(i),\ i=1,...,k_n$, where $\xi_n(i) = \{T^jA_n(i)\ |\
j=0,...,h_i-1\}$. Given $\e > 0$ choose $n$ sufficiently large so that
${\rm diam}(A_n) <\e/2$ and ${\rm diam}(T^{-1}A_n) <\e/2$. For  $n$,
consider two collections of clopen subsets:
$(T^{h_1-1}A_n(1),...,T^{h_{k_n}-1}A_n(k_n))$ and $(A_n(1),...,A_n(k_n))$.
Let $S_i$ be a one-to-one continuous map from $T^{h_i-1}A_n(i)$ onto
$A_n(i+1),\ i=1,...,k_n-1$. Let $Y_1 =\Om \setminus T^{h_{k_n}-1}A_n(k_n)$
and $Y_2 =\Om \setminus A_n(1)$. Define a one-to-one continuous map $S$
from $Y_1$ onto $Y_2$ by
$$
Sx =\left\{\begin{array}{ll} Tx, & {\rm if}\ x\notin
\bigcup_{i=1}^{k_n}T^{h_i-1}A_n(i)\\ S_ix,&{\rm if}\ x\in T^{h_i-1}A_n(i),\
i=1,...,k_n-1
\end{array} \right.
$$
We note that $d(Tx,Sx) \le {\rm diam}(T^{-1}A_n),\ x\in Y_1$ and $d(T^{-1}x,
S^{-1}x) \le {\rm diam}(A_n),\ x\in Y_2$.

This construction defines an $S$-tower with base $A_n(1)$ and  top
$T^{h_{k_n}-1}A_n(k_n)$. Next, by using  cutting and stacking,
we may extend $S$ to a homeomorphism of $\Om $ denoted also by $S$.
Clearly, $S$ has rank 1. By our construction,
$$
d_w(S,T)=\sup_{x \in \Om}d(Sx, Tx)+\sup_{x\in \Om}d(S^{-1}x,T^{-1}x) <\e.
$$
This proves that $T$ can be approximated in $\tau_w$ by a homeomorphism $S$
of rank 1. Thus we have shown that $\Min \subset
\overline{\R(1)}^{\tau_w}$.

Now we will prove that $\Min \subset \overline{\R(1)}^\tau$. Let
$\mu_1,...,\mu_m$ be a given collection of Borel measures and let $\e >0$.
Take $T\in \Min$ and construct $(A_n)$ and $(\xi_n)$ as in the first part
of the proof. It was proved in Theorem \ref{per2} that $n$ can be
chosen so large that $\mu_i(A_n) <\e/2$ and $\mu_i(T^{-1}A_n) <\e/2, \
i=1,...,m$. We apply the above method to define a
homeomorphism $S$ of $\Om$ of rank 1. Obviously, $E(S,T) \subset A_n\cup
T^{-1}A_n$. Then $\mu_i(E(S,T)) < \e $ for $i =1,...,m$. Therefore $S\in
U(T;\mu_1,...,\mu_m;\e)\cap \R(1)$.\hfill$\square$ \\

The next statement can be proved by using the techniques of proof of
Theorem \ref{min}.

\begin{corollary}\label{fin-min}  Let $T$ be a homeomorphism of $\Om$
which has a finite decomposition into minimal components. Then $T\in
\overline{\R(1)}^\tau$ and $T\in \overline{\R(1)}^{\tau_w}$.
\end{corollary}

It follows from Theorems \ref{r1-min}  and \ref{min}  that the following
result holds.

\begin{corollary}\label{all} $(1)$  $\od$ is a
dense $G_\de$-set in $ \overline{\Min}^{\tau_w}$, i.e. $\od$ is a set of
the second category in $\overline{\Min}^{\tau_w}$.\\ $(2)$ $\od$ is a dense
$F_{\sigma\de}$-set  in $\overline{\Min}^\tau$.\\
$(3)$ $\Min$ is closed neither in $\tau$ nor in $\tau_w$.
\end{corollary}

It follows from Corollary \ref{all}  that a typical minimal homeomorphism
is saturated (see Remark \ref{sat}).

\begin{theorem}\label{min-mix1} $\mix\cap \overline{\Min}^{\tau_w}$  is
nowhere dense in $(\overline{\Min}^{\tau_w}, \tau_w)$.
\end{theorem}

\noindent{\it Proof}. Let $T\in \overline{\Min}^{\tau_w}$. Take a
$\tau_w$-neighborhood $W(T;F_1,...,F_m)$. We showed in Theorem \ref{min}
that there exists an odometer $S \in W(T;F_1,...,F_m)$. Let $\F =
(E,SE,...,S^{n-1}E)$ be an $S$-tower such that $\bigcup_{i=0}^{n-1} S^iE
=\Om$ and $S^nE = E$.  Then  $W(S; E,SE,...,S^{n-1}E)$ consists of the
homeomorphisms $R\in\h$ such that $R(S^iE) =S^{i+1}E,\ i=0,1,...,n-1$.
Note that, choosing $n$ sufficiently large, we have that every $F_i$ is
an $\F$-set. Therefore $W(S; E,SE,...,S^{n-1}E)  \subset W(T;F_1,...,F_m)$
since $SF_i = TF_i$ and each $F_i$ is a union of some atoms from $\F$.
Then for any $R W(S; E,SE,...,S^{n-1}E)$, we have  $\{m\in \Z\ |\
R^mE\cap E \neq \emptyset\} = n\Z$ and therefore $R$ cannot be mixing. This
proves that $\mix$ is nowhere dense in $\overline{\Min}^{\tau_w}$.
\hfill$\square$
\\

\noindent{\bf Remark} Since $\mix\cap\od =\emptyset$, Theorem \ref{min-mix1}
is consistent with the conclusion of Corollary \ref{all} that $\od$ is a set
of second category in $\overline{\Min}^{\tau_w}$.
\\

We now introduce a concept which will allow us to characterize the closure
of $\Min$ in the weak topology $\tau_w$.

\begin{definition}\label{mov-def}
We call a homeomorphism $T$ moving if for any non-trivial clopen set $F$
each of the sets $TF\setminus F$ and $F\setminus TF$ is not empty. A
homeomorphism $T$ is called weakly moving if $TF\neq F$ for every
non-trivial $F\in CO(\Om)$.
\end{definition}

We denote by $\mov$ and w-$\mov$ the sets of moving and weakly
moving homeomorphisms, respectively. Obviously, $\mov\subset$
w-$\mov$ and if $T\in \mov$ (or $T\in $ w-$\mov$) then also
$T^{-1} \in \mov$ (or $T^{-1}\in $ w-$\mov$).

It follows immediately that $\Min$ and $\mix$ are subsets of
w-$\mov$. We make this statement more precise in Theorem \ref{mov} below.

Note that the set w-$\mov$ is closed in $\tau_w$. Indeed, it is easily seen
that
\begin{equation}\label{wmov}
w\text{-}\mov = \bigcap_{F\in CO(\Om)} W({\mathbb I};F)^c.
\end{equation}
Another simple observation from (\ref{wmov}) says that if $P\in \per_0$,
then $P\notin$ w-$\mov$.

\begin{theorem} \label{mov} $T\in \mov\ \iff\ T\in
\overline{\od}^{\tau_w}\ =\ \overline{\Min}^{\tau_w}$.
\end{theorem}

\noindent {\it Proof}. Let $T\in \overline{\od}^{\tau_w}$.
Consider a neighborhood $W(T;F)$ where $F\in CO(\Om)$. Then
$W(T;F)$ contains an odometer $S$ such that $SF=TF$. Let $\xi$ be
an $S$-tower such that $F$ is a $\xi$-set. Since every atom of
$\xi$ is shifted by $S$ to the next level of the tower, we get
that $TF\setminus F = SF \setminus F\neq \emptyset$ and
$F\setminus TF = F \setminus SF\neq \emptyset$.

Conversely, suppose that $T$ is moving. Let $W(T; F_1,...,F_n)$ be
a neighborhood where $\zeta= (F_1,...,F_n)$ is a partition of $\Om$ into
clopen sets. We know that for every $i$ the sets $TF_i\setminus
F_i$ and $F_i\setminus TF_i$ are non-empty. Take the intersection
$\zeta\wedge T(\zeta)$ of partitions $\zeta$ and $T\zeta$ consisting of
atoms $TF_i\cap F_j= F_{ij}$ where $i,j= 1,...,n$ and some of $F_{ij}$'s
may be empty.

We first consider a particular case when for every $i=1,...,n$,
\begin{equation}\label{card}
|\{1\le j\le n : TF_i\cap F_j\neq \emptyset\}| = |\{1\le j\le n : TF_j
\cap F_i\neq \emptyset\}|
\end{equation}
Let $A= (a_{ij})$ be an $n\times n$ matrix where $a_{ij} =1 $ if
$TF_i\cap F_j\neq\emptyset$ and $a_{ij} =0$ otherwise. In other
words, we consider the directed graph $\G$ with the set of
vertices $(1,...,n)$ and the set of arrows defined by $A$: an arrow  from
$i$ to $j$ exists if and only if $a_{ij} =1$. Note that $\G$ may have
loops, i.e. arrows that begin and end at the same vertex. In other words,
relation  (\ref{card}) says that the number of arrows coming to a vertex
$i$ equals the number of arrows outgoing from $i$. We claim that
$(\G,A)$ is a connected graph. Indeed, let $Z(i)$ be the set of vertices
that can be connected with $i$ by a path. We show that $Z(i)= (1,...,n)$
for every $i$. If $j\in Z(i)$ then there exist $j_0 = i, j_1,...,j_s = j$
such that $TF_{j_k} \cap F_{j_{k+1}}\neq \emptyset,\ k = 0,...,s-1$. Assume
that $\La = (1,...,n)\setminus Z(i)\neq \emptyset$ and denote by $E
=\bigcup_{j\in \La}F_j$. Since $T\in \mov$, we get that $TE^c \setminus E^c
\neq \emptyset$ or $TE^c\cap E \neq\emptyset$ which contradicts the
assumption.

Since we have a connected graph $(\G,A)$ satisfying (\ref{card}), we
can choose an Euler path $L$, consisting of arrows, which goes through all
vertices and takes every arrow only once (see [G] for details).

We now construct a homeomorphism $S$ of $\Om$ using the path $L$. To every
arrow from $i$ to $j$ we associate the set $F_{ij}$. We start with some
vertex $i_0$ and let $S$ be a homeomorphism from $F_{i_0j}$ onto $F_{jk}$ if
the arrow from $j$ to $k$ follows that from $i_0$ to $j$ in $L$. Then we
extend the definition of $S$ going along $L$. Thus, $S$ is defined on $\Om$
and atoms of $T(\zeta)\wedge\zeta$ form an $S$-tower, and since the path $L$
is annular, the top of the tower is mapped by $S$ onto the base. Moreover,
one can verify that by definition of $S$, $TF_i = SF_i$ for all $i$. Indeed,
fix some $i$ and consider the sets $F_{ji}$ and $F_{ik}$ when $j,k$ run over
$(1,...,n)$. We see that $\bigcup_jF_{ji} = F_i$ and $\bigcup_kF_{ik}=
TF_i$. Hence $S$ maps $F_i$ onto $TF_i$. Therefore we have found an odometer
$S$ which belongs to the neighborhood $W(T; F_1,...,F_n)$.

In general, (\ref{card}) does not hold. But we can slightly modify the above
construction to obtain the result. We note that there are positive integers
$m_{ij}$ such that
\begin{equation}\label{m}
\sum_{\{j :\ a_{ij} =1\}} m_{ij} = \sum_{\{j :\ a_{ji} =1\}}
m_{ji}.
\end{equation}
This means that we can consider a new graph $(\G, A')$ over the same set of
vertices $(1,...,n)$ and with an extended set of arrows: if $i$ and $j$ are
such that $a_{ij} =1$, then we take exactly $m_{ij}$ directed arrows from
$i$ to $j$. Given two connected vertices $i$ and $j$, we can assign a number
from 1 to $m_{ij}$ to each arrow from $i$ to $j$. In this way, we see by
(\ref{m}) that the connected graph $(\G,A')$ has the following property: the
number of arrows arriving at each vertex is the same as the number of arrows
leaving that vertex. Therefore, there exists a closed path $L'$ that goes
through all vertices (visiting each vertex several times) and includes each
arrow only once.

To construct an odometer $S$, we divide every non-empty set
$F_{ij}$ into $m_{ij}$ non-empty clopen subsets $E_{ij,k},\
k=1,...,m_{ij}$. To define  $S$, we start with a vertex $i_0$
and some outgoing arrow $l(i_0,i_1;k_1)$ from $i_0$ to $i_1$ with the
assigned number $1\le k_1 \le m_{i_0i_1}$. Then $S$ is defined by
the following rule: if the arrow $l(i_1,i_2,k_2)$ from $i_1$ to $i_2$ with
the number $k_2$ is next to $l(i_0,i_1;k_1)$ in $L'$ the we set  $S :
E_{i_0i_1,k_1} \to E_{i_1i_2,k_2}$. One can continue this procedure going
along  $L'$ until the last arrow in $L'$ has been used. This arrow has the
largest number $m_{i_{n-1}i_0}$ amongst those that return to $i_0$.
As before, it is easy to check that the sets $(E_{ij,k}\ i,j
=1,...,n;k= 1,...,m_{ij})$ define an $S$-tower and $SF_i =TF_i$ for all
$i$. \hfill$\square$
\\

Let ${\cal T}t$ denote the set of all topologically transitive
homeomorphisms of $\Om$.

\begin{corollary}\label{tt} $(1)$ ${\cal T}t \subset \mov$ and
$\overline{{\cal T}t}^{\tau_w} = \overline{\Min}^{\tau_w}$ ;\\
$(2)$ $\overline{\Min}^\tau \neq\overline{\Min}^{\tau_w}$.
\end{corollary}

\noindent{\it Proof}. It is easily seen that every $T\in {\cal T}t$
satisfies Definition \ref{mov-def}. Thus, we have
$$
\mov\ \supset\ {\cal T}t \ \supset\ \Min,
$$
and the result follows from Theorem \ref{mov}.

On the other hand, if  some $T$ from ${\cal T}t$ has a periodic orbit then
such a $T$ cannot be in $\overline{\Min}^\tau$ in view of Theorem
\ref{ap}.\hfill{$\square$}

\begin{corollary} $\Min \subset \overline{\od}^{\tau_w} =
\overline{\Min}^{\tau_w} \subset$ w-$\mov$ and  $\mix \subset
\overline{\od}^{\tau_w} = \overline{\Min}^{\tau_w}  \subset $ w-$\mov$.
\end{corollary}

\noindent{\it Proof}. We need only show that if $T$ is either minimal or
mixing then $T$ belongs to $\overline{\od}^{\tau_w}$. For this, assume that
it is not true, i.e. $T\notin \overline{\od}^{\tau_w} = \mov$. Then there
exists a proper clopen subset $F$ such that either (i) $TF=F$ or (ii) $TF
\subset F$ or (iii) $F\subset TF$. If (i) holds then $T$ is neither minimal
nor mixing. If (ii) holds, then $T^{n+1}F \subset T^nF, \ n \in \N $, and
the set $\bigcap_{n\ge 0} T^nF$ is a closed $T$-invariant set. Thus, $T$
cannot be minimal. Denote by $E= F\setminus TF$. Then $E, TE, ...,
T^nE,...$ are pairwise disjoint and therefore $T$ cannot be mixing. Similar
arguments are used in case (iii).\hfill$\square$
\\

Analyzing the proof of Theorem \ref{mov} we can immediately deduce the
following consequence.

\begin{theorem}\label{per0-1} A homeomorphism $T$ of $\Om$ belongs to
$\overline{\per_0}^{\tau_w}$ if and only if for each clopen $F$ either
$TF=F$ or $TF\setminus F \neq\emptyset$ or  $F\setminus TF \neq\emptyset$.
Therefore every homeomorphism $T\in
\overline{\Min}^{\tau_w}$ can be approximated by a periodic
homeomorphism $P\in \per_0$ in the topology $\tau_w$, that is
$\overline{\Min}^{\tau_w}\subset \overline{\per_0}^{\tau_w}$.
\end{theorem}

\noindent {\it Proof}. From the proof of Theorem \ref{mov}, we see that if
the hypotheses of the Theorem hold, one can construct a periodic
homeomorphism $S$ of $\Om$ using the Euler path $L$ in the same way as for
odometers. \hfill$\square$
\\

A simple consequence of Theorem \ref{per0-1} is the following fact:
$\overline{\per}^{\tau_w} = \overline{\per_0}^{\tau_w}$.

\begin{remark}\label{rem} {\rm We may also use Theorem \ref{mov} to show
that for a minimal homeomorphism $T\in \h$, $[[T]]_{ap} \subset
\overline{\Min}^\tau$ as stated in Theorem \ref{ap-gamma}. For this, we need
to prove that if $\g \in [[T]]_{ap}$, then $\g \in \mov$, that is for any
clopen $E$ the sets $\g E\setminus E$, $E\setminus \g E$ are non-empty.
Assume that there exists $F \in CO(\Om)$ such that $\g F \subseteq F$ and
deduce from this that $\g$ must have a periodic part. Choose a
Kakutani-Rokhlin partition $\xi$ such that $\g\in \G(\xi)$ and the clopen
sets $F$ and $\g F$ are $\xi$-sets (we use here notation of Section 2). Then
$F$ and $\g F$ are unions of some atoms $D^j(v)$ of $\xi$. We have two
possibilities: either every $D^j(v)\subset F$ does not leave the subtower
$\xi(v)$ under the action of $\g$, or there exists $D^{j_0}(v_0)$ in $F$
such that $D^{j_0}(v_0)$ belongs to some $L_{ti}(J)$ (or $L_{bk}(J'))$. From
the first case we get that $\g$ must have a periodic part inside of
$\xi(v)$. The second case implies that if $D^{j_0}(v_0)\subset F \cap
L_{ti}(J)$, then $L_{ti}(J)\subset F$ (as well as $\g(L_{ti}(J))$) is also a
subset in $F$ since the set $\g F$ is a $\xi$-set. Then $\g$ again has a
periodic part. }
\end{remark}

We can strengthen the first statement of Theorem \ref{min} by giving a
complete description of closures of various classes of homeomorphisms in the
topology $\tau$.

\begin{theorem}\label{aper} $ \overline{\od}^\tau = \overline{\R(1)}^\tau =
\overline{\Min}^\tau = \overline{{\cal T}t}^\tau \cap \ap=\ap$.
\end{theorem}
{\it Proof}. Let $T$ be an aperiodic homeomorphism of $\Om$ and let $U(T) =
U(T;\mu_1,...,\mu_k;\e)$ be a $\tau$-neighborhood of $T$.

We will apply the following result established in [BDM] to prove the Rokhlin
lemma.

\begin{lemma}\label{cover} Let $T$ be an aperiodic homeomorphism of a Cantor
set $\Om$, $\mu_1,...,\mu_k\in \M_1(\Om),\ \e>0$. Given a positive integer
$n\ge 2$, there exists a partition of $\Om$ into a finite number of clopen
$T$-towers $(\eta_1,\ldots,\eta_q)$ such that the height $h(\eta_i)$ of
every tower is at least $n$. Moreover, these towers can be chosen such that
\begin{equation}\label{rokh}
\mu_i(\bigcup_{j=0}^{n-1}T^{-j}B)> 1 -\e,
\end{equation}
where $B= \bigcup_{i=1}^q B_i$ and $B_i$ is the base of $\eta_i$.
\end{lemma}
{\it Sketch of proof}. We begin with a clopen finite disjoint cover
$(U_1,...,U_k)$ of $\Om$ such that $T^jU_i \cap U_i =\emptyset,\
j=1,...,n-1,\ i=1,...,k$. Consider $\xi_1 = (U_1, TU_1,...,T^{n-1} U_1)$ and
define $C_1 = \bigcup_{j=0}^{n-1}T^jU_1$. Let $U_i^1 = U_1 - C_1,\
i=2,...,k$. Define
$$
U_2^1(i)=\{x\in U_2^1\,:\, T^ix\in U_1,\ T^jx\notin U_1,\ 0\leq j\leq
i-1\},\ i=1,\ldots,n-1,
$$
and
$$
U_2^1(0)=\{x\in U_2^1\,:\,T^jx\notin U_1,\mbox{ for all }1\leq j\leq n-1\}.
$$
Each set $U_2^1(i)$ is the base of the $T$-tower
$$
\xi_2^1(i)=\{U_2^1(i),TU_2^1(i),\ldots,T^{n-1+i}U_2^1(i)\}, \mbox{ for all
}i=0,1,\ldots,n-1.
$$
Take the set $U_1^1 = U_1 \setminus\ \bigcup\limits_{i=1}^{n-1}T^iU_2^1(i)$
as the base of a subtower $\xi_1^1$ of $\xi_1$. We get the collection of
disjoint $T$-towers $\Xi(1) =
\{\xi_1^1,\xi_2^1(0),\xi_2^1(1),\ldots,\xi_2^1(n-1)\}$ each of which is of
height at least $n$. Denote by $C_1^1, C_2^1(0), C_2^1(1),...,C_2^1(n-1)$
the supports of the corresponding towers from $\Xi(1)$. Note that $C_1^1\cup
C_2^1(0)\cup C_2^1(1)\cup\cdots \cup C_2^1(n-1) = C_1 \cup C_2^1$ where
$C_2^1=\bigcup\limits_{s=0}^{n-1}T^i U_2^1$. Define the sets
$U_i^2=U_i^1\setminus\ C_2^1,$ for all $i=3,4,\ldots,k$. For each tower
$\xi$ from $\Xi(1)$, denote by $U_3^2(\xi)$ the subset of $U_3^2$ which
consists of those points whose $T$-orbits meet $\xi$ for at most $n-1$
iterations. We can now apply the construction used above to the set
$U_3^2(\xi)$ and the tower $\xi\in\Xi(1)$. Repeating this procedure at most
$k-1$ times we will finally obtain a collection $\Xi =(\xi_1,...,\xi_s)$ of
disjoint $T$-towers which covers $\Om$ such that the height of each tower
$\xi\in \Xi$ is at least $n$.

Choose $m\in \N$ such that $1/m < \e$ and let $V_j$ denote the top of
$\xi_j$. Using the obvious fact that amongst $m$ pairwise disjoint subsets
of $\Om$ at least one must have $\mu_i$-measure less than $\e$, we may
choose a clopen set $B= T^{-K} V$, where $0\le K< n$ and $V= \bigcup_{j=1}^s
V_j$, such that (\ref{rokh}) holds. To get the $T$-towers
$\eta_1,...,\eta_q$, one refines the existing towers $\xi_1,...,\xi_s$ such
that the base of each $\eta_i$ is a subset of $B$ and the top is a subset of
$T^{-1}B$. Full details of the proof are in [BDM]. \hfill$\square$
\\

We now return to the proof of the Theorem. Define a homeomorphism $S$ as
follows. Let $S_i$ be a homeomorphism sending the top $T^{h(\xi_i) -1}B_i$
of $\xi_i$ onto the base $B_{i+1}$ of $\xi_{i+1}$ ($i=1,...,m-1)$ and $S_n$
maps $T^{h(\xi_m) -1}B_m$ onto $B_1$. Take $n=2$ in Lemma \ref{cover} and
define
$$
Sx=\left\{
\begin{array}{ll}
Tx, & \mbox{ if }x\in
\bigcup\limits_{i=1}^m\bigcup\limits_{j=0}^{h(\xi_i)-2}T^jB_j,
\\
S_ix & \mbox{ if }x\in T^{h(\xi_i)-1}B_i\mbox{ for some } i= 1,...,m.
\end{array}\right.
$$
It follows from (\ref{rokh}) that the homeomorphism $S$ belongs to $U(T)$.
On the other hand, there exists a homeomorphism $S_0\in \od$ such that $S_0x
\neq Sx$ only if $x\in T^{h(\xi_m) -1}B_m$, hence $S_0\in U(T)$. This proves
that $\overline{\od}^\tau = \ap$. \hfill$\square$.
\\

{\it Acknowledgement}. We would like to thank A.~Kechris and G.~Hjorth for
helpful discussions of our results and especially for their contributions to
the proof of Theorem \ref{property}. The first named author thanks the
University of New South Wales for the warm hospitality and the Australian
Research Council for its support. S.B. and J.K. acknowledge also the support
of the Torun University and the Kharkov Institute for Low Temperature
Physics during exchange visits.

\vskip1cm \noindent {\small {\it S.~Bezuglyi\\
Department of Mathematics\\
Institute for Low Temperature Physics\\
47, Lenin ave.\\
61103 Kharkov, UKRAINE\\
\smallskip
\noindent bezuglyi@ilt.kharkov.ua}}
\\
\\
\noindent {\small {\it A.H.~Dooley\\
School of Mathematics\\
University of New South Wales\\
Sydney, 2052 NSW, AUSTRALIA\\
\smallskip
\noindent tony@maths.unsw.edu.au}}
\\
\\
\noindent {\small {\it J.~Kwiatkowski\\
Faculty of Mathematics and Computer Sciences\\
Nicholas Copernicus University\\
ul. Chopina 12/18\\
87-100 Torun, POLAND\\
\smallskip
\noindent jkwiat@mat.uni.torun.pl}}

\end{document}